%% file: oco_perturbed.tex
\documentclass[11pt]{article}
\usepackage{geometry}                
\usepackage[utf8]{inputenc} 
\usepackage[T1]{fontenc}    
\usepackage{url}            
\usepackage{booktabs}       
\usepackage{amsfonts}       
\usepackage{nicefrac}       
\usepackage{microtype}      
\usepackage{graphicx}
\usepackage{algorithm}
\usepackage{algorithmic}
\usepackage{wrapfig}
\usepackage{fullpage}
\usepackage{authblk}

\title{Online Convex Optimization with Perturbed Constraints}

\author[$*$]{V\'ictor Valls}
\author[$\dagger$]{George Iosifidis}
\author[$\dagger$]{Douglas J.\ Leith}
\author[$*$]{Leandros Tassiulas}
\affil[$*$]{Department of Electrical Engineering and Institute for Network Science, Yale University}
\affil[$\dagger$]{School of Computer Science and Statistics, Trinity College Dublin}

\date{}

\usepackage{mathtools}
\usepackage{amsthm}
\usepackage{dsfont}
\usepackage{xcolor}
\graphicspath{{Figures/}}

\def\N{{\rm I\!N}}
\def\R{{\rm I\!R}}
\def\B{B_\psi}

\def\1{\mathds{1}}

\newtheorem{lemma}{Lemma}
\newtheorem{theorem}{Theorem}
\newtheorem{corollary}{Corollary}

\newtheorem{proposition}{Proposition}

\begin{document}

\maketitle 

\begin{abstract}
This paper addresses Online Convex Optimization (OCO) problems where the constraints have additive perturbations that (i) vary over time and (ii) are not known at the time to make a decision. Perturbations may not be i.i.d.\ generated and can be used to model a time-varying budget or commodity in resource allocation problems. The problem is to design a policy that obtains sublinear regret while ensuring that the constraints are satisfied on average. To solve this problem, we present a primal-dual proximal gradient algorithm that has $O(T^\epsilon \vee T^{1-\epsilon})$ regret and $O(T^\epsilon)$ constraint violation, where $\epsilon \in [0,1)$ is a parameter in the learning rate. Our results match the bounds of previous work on OCO with time-varying constraints when $\epsilon = 1/2$; however, we (i) define the regret using a time-varying set of best fixed decisions; (ii) can balance between regret and constraint violation; and (iii) use an adaptive learning rate that allows us to run the algorithm for any time horizon.
\end{abstract}

\section{Introduction}


The Online Convex Optimization (OCO) framework was introduced in \cite{Zin03} and it is widely used to model applications such as spam filtering, portfolio selection, recommendation systems, among many others \cite{Haz16}. In short, OCO consists of a sequence of games where in each round $t \in \N$ an agent selects an action $x_t$ from a convex set $X \subset \R^n$ and suffers a cost $f_t(x_t)$, where $f_t : \R^n \to \R$ is convex. 
Crucially, the cost function is not known at the time of making a decision, and it may even be selected by an adversary after the action has been played. The goal is to design a policy or algorithm that selects a sequence of actions $\{x_t\}$, $t =1 ,\dots,T$ from $X$ so that the regret
\begin{align}
R(T) & := \sum_{t=1}^T f_t(x_t) - \min_{x \in X} \sum_{t=1}^T f_t(x)
\end{align}
increases sublinearly, i.e., ${\lim \sup}_{T \to \infty} R(T) /T \le 0$. Hence, the incurred cost is asymptotically as good as the best fixed decision in hindsight.\footnote{The regret captures the difference between the incurred cost and the cost obtained by an ``offline'' algorithm that has knowledge of all the cost functions from $t=1,\dots,T$. The decision of the offline algorithm is, however, more restricted, as it can only choose one vector from $X$. }

\cite{Zin03} showed that the online gradient descent (OGD) algorithm can obtain sublinear regret when the action set $X$ is bounded and the convex cost functions $f_t$, $t\in \N$ have bounded subgradients. 
The algorithm consists of update
\begin{align}
x_{t+1} 
& = \mathcal P_X ( x_t  - \alpha_t f'_t(x_t) ), \qquad t = 1,2,\dots
\label{eq:ogd}
\end{align}
where $\alpha_t = 1 / \sqrt t$ is the learning rate (or, step size), $f_t'(x_t)$ a subgradient of the previous cost function at $x_t$, and $\mathcal P_X$ the Euclidean projection onto the convex set $X$. 
Note from Eq.\ \eqref{eq:ogd} that action $x_{t+1}$ is selected using only information available at time $t$. 


\subsection{OCO with long-term (LT) constraints}
\label{sec:avgconst}

In the standard OCO setting, set $X$ encompasses \emph{all} the constraints that an online policy and a fixed decision in hindsight must satisfy. However, sometimes it is useful to treat constraints differently depending on whether they are (i) \emph{instantaneous constraints} that have to be satisfied in every iteration, or (ii) \emph{long-term constraints} that must be satisfied only in the long-term or average sense (the formal definition is given in Sec.\ \ref{sec:model}).

The two most prominent reasons for considering long-term constraints are the following. First, \underline{\smash{more freedom of choice}}. By not requiring that the long-term constraints are satisfied in every iteration, it is possible to devise  policies with specific properties (e.g., lower-complexity per iteration; see Sec.\ \ref{sec:relatedwork}) and to model new classes of OCO problems. 
%
For instance, in wireless communications systems, the power needed to transmit a message is not known a priori (as it depends on the channel conditions, the behavior of other users, etc.) and the goal is to adjust the transmission power to maximize the rate while keeping the \emph{average} power consumption below a predefined threshold \cite{MTY09}. That is, it is possible to transmit a message using more power than what it is allowed to use on average as long as the average power consumption constraint is met. 

The second reason is that they allow us to handle \underline{\smash{online constraints}}. 
More specifically, constraints that (i) change over time and (ii) are not known by the decision maker at the time of making a decision.
For example, in online network flow problems with ``offline'' constraints (e.g., online shortest path routing \cite[pp.\ 7]{Haz16}) the amount of flow $x_t$ allocated to each of the links has to satisfy the equation $Ax_t = b$ in each time $t\in \N$,\footnote{See \cite{Roc84,GNT06} for an introduction to modeling different types of network flow problems.} where $A$ is the routing matrix and $b$ a request vector that indicates the supply/demand of flow (of material, traffic, information, etc.) at each of the nodes. 
With online constraints, we have $b_t$ instead of $b$ (i.e., $Ax_t = b_t$) and must select $x_t$ {without} knowledge of $b_t$. Hence, it is not possible to design an algorithm that guarantees that $Ax_t$ is equal to  $b_t$ in every iteration.
An example of this type of constraint is in power allocation problems in data centers where the number of machines that run at a given time  ($x_t$) has to be decided before the real workload ($b_t$) can be observed \cite{GIN19}.\footnote{See \cite[Sec.\ 8]{MTY09} for a similar example with CPUs.}
It is important to emphasize that like the cost functions in standard OCO, the perturbation can be a function of the past actions. For instance, in network flow problems \cite{GNT06}, the previous resource allocation decisions ($x_{t-1} , \dots, x_1$) affect the quality of experience of the users and therefore, how those generate resource requests $(b_t)$.
Another example is online display advertising \cite{AD15} where the perturbations $ b_t $ represent a budget that varies with time and depends on the past rewards. Or more precisely, past actions affect the rewards and, therefore, the future advertising budget.

\subsection{Contributions and related work}
\label{sec:relatedwork}

In this paper, we consider OCO problems where the long-term constraints have additive perturbations $b_t$ that (i) change over time and (ii) are not known by the decision maker at the time of making a decision $x_t$.
%
%
%
We do not require the perturbations to be i.i.d.\  or to have any other statistical property (see assumptions in Sec.\ \ref{sec:assumptions}). 
The problem is to design and algorithm that obtains sublinear regret and ensures that the constraints are satisfied on average. 
The problem addressed is important because non-i.i.d.\ perturbations can be used to model more accurately a time-varying budget or commodity in online resource allocation problems; for instance, when the actions that the agent made in the past affect the future constraints (i.e., the perturbations).

\begin{wrapfigure}{r}{0.3\textwidth}
 \small \begin{center}
    \includegraphics[width=0.2\textwidth]{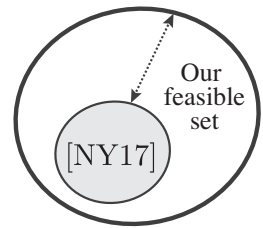}
  \end{center}
  \caption{ {Showing how our feasible set compares to the feasible set in   \cite{NY17} (shaded area). 
  }}
  \label{fig:sets}
\end{wrapfigure}

We solve the OCO problem with perturbed constraints using a novel primal-dual proximal gradient algorithm (Algorithm \ref{al:roco}) that obtains $O(T^\epsilon \vee T^{1-\epsilon})$ regret and $O(T^\epsilon)$ accumulated constraint violation, where $\epsilon \in [0,1)$ is a parameter in the learning rate. The best decision in hindsight is selected from a feasible set\footnote{i.e., a set where the best fixed decisions satisfies the constraints on average.} that changes over time depending on the perturbations. 
Our algorithm allows us to balance between regret and constraint violation by choosing $\epsilon$ accordingly. 
When $\epsilon = 1/2$, our bounds match the best well-known rates \cite{NY17}, but we can also obtain faster violation rate than $O(\sqrt T)$. For example, with $\epsilon = 1/4$, we have $O(T^{3/4})$ regret and $O(T^{1/4})$ constraint violation. Another key characteristic of our algorithm is that the learning rate is adaptive. Hence, we do not need to fix in advance the time the algorithm will run (which is not known in many resource allocation problems). Furthermore, adaptive learning rates are preferable than extending the horizon with the ``doubling trick'' \cite[Sec.\ 2.3]{CL06}.
Table \ref{table:compare} shows a summary of how our technical approach compares to previous approaches. However, we must emphasize that the problem we address in this paper is fundamentally different from the one considered in previous work; especially, \cite{MJY12,JHA16}. Our motivation for relaxing the constraints is not the complexity of the projection onto the feasible set, but actually that the feasible set is not known and changes over time (because of the perturbations). 

\begin{table}
\label{table:compare}
{{
\caption{The \emph{feasible set} is the set from which the fixed decision in hindsight is selected. In \cite{JHA16}, $\beta \in (0,1)$. In \cite{YNW17}, the bounds hold with probability at least $1 - \delta$ where $\delta \in (0,1)$, and the feasible set is defined in expectation. In this work, $\epsilon \in [0,1)$. The works with ($\dagger$) consider only static (or ``offline'') constraints.}
}}
\centering
{\small
\begin{tabular}{l  cccc}
\toprule
\textbf{Paper}& \textbf{Feasible Set}  & \textbf{Regret}& \textbf{Constraint Violation} & \textbf{Learning Rate} \\
\hline
\cite{MJY12}$^\dagger$  &  fixed  &  $O(T^{\{\frac{1}{2}, \frac{2}{3}\}}) $ & $O(T^{\{\frac{3}{4}, \frac{2}{3}\}}) $ & constant   \\
\cite{JHA16}$^\dagger$  &  fixed  & $O(T^{\beta} \vee T^{1-\beta})$ & $O(T^{1-\beta/2})$ & adaptive   \\
\cite{NY17} & fixed & $O(\sqrt T)$ & $O(\sqrt T)$ & constant \\
\cite{YNW17}& \text{ \ \ }fixed$^*$  \   & $O(\sqrt T \log(T) \log(\frac{1}{\delta})) $ & $O(\sqrt T \log(T) \log^{3/2}(\frac{1}{\delta}))$ & constant \\
\hline
This work &  time-varying &  $O(T^\epsilon \vee T^{1-\epsilon})$ & $O(T^\epsilon)$ & adaptive \\
\bottomrule
\end{tabular}
}
\end{table}


\textbf{Related work:}
The first works of OCO with long-term constraints were motivated by the complexity of the projection step in OGD. In brief, when set $X$ is composed of general convex constraints, the projection step involves solving a convex program that can be computationally burdensome. For example, projections onto the semidefinite cone. 
Expensive projections are dealt with in offline convex optimization by carrying them out only in the last iteration \cite{MYJ+12} or less often \cite{CGP16}. However, such approaches cannot be used in OCO problems as every action played incurs an instantaneous cost.
The latter was noted in \cite{MJY12}, which formalized the OCO problem with static long-term constraints and proposed two algorithms based on variational inequalities \cite{Nem94}. First, a gradient-based algorithm that obtains $O(\sqrt T)$ regret and $O(T^{3/4})$ constraint violation for general OCO problems; and second, a mirror-prox algorithm that obtains $O(T^{2/3})$ regret and constraint violation when the constraints are polyhedral. 
The paper in \cite{JHA16} extends the work \cite{MJY12} by proposing an algorithm for general OCO problems with long-term constraints that can balance regret and constraint violation. In particular, the algorithm obtains  $O(T^{\beta} \vee T^{1-\beta})$ regret and $O(T^{1-\beta/2})$ constraint violation where $\beta \in (0,1)$ is a design parameter.
Furthermore,  and unlike in \cite{MJY12}, the learning rate is adaptive, and so the algorithm can run for any time horizon.

Regarding online constraints, the work in \cite{MTY09} considers online learning problems with constraints that can vary in an arbitrary and possibly adversarial manner. The paper shows that the highest reward-in-hindsight while satisfying the constraints is not attainable in general, except for the case where the feasible set (i.e., the set from which the best fixed decision in hindsight is selected) is convex. 
The latter result motivated the work in \cite{NY17} to consider OCO problems with convex time-varying constraints (not only with additive perturbations) and proposed an algorithm that obtains $O(\sqrt T)$ regret and constraint violation.
However, the performance of the proposed algorithm is compared to the best fixed decision in hindsight that satisfies \emph{every} time-varying constraint. This is in stark contrast to our work where the feasible set changes over time depending on the perturbations. Furthermore, our feasible set always contains the feasible set in \cite{NY17} (see Fig.\ \ref{fig:sets} and Sec.\ \ref{sec:model}), which means that we compare the cost of our algorithm with a larger set of (possibly better) fixed decisions in hindsight. Hence, the cost of the best fixed decision in hindsight in our work can be smaller than in \cite{NY17}.
%
%
Finally,  \cite{YNW17} considers online constraints that are i.i.d.\ generated where the feasible set is defined in expectation. The proposed algorithm obtains $O(\sqrt T)$ regret and constraint violation in expectation, and $O(\sqrt T \log(T) \log(\frac{1}{\delta}))$ regret and $O(\sqrt T \log(T) \log^{3/2}(\frac{1}{\delta}))$ constraint violation bounds that hold for every sample path with probability $1 - \delta$, $\delta \in (0,1)$. 
The algorithms in \cite{NY17,YNW17} are based on Lyapunov techniques from stochastic network optimization \cite{Mey08, GNT06}; the learning rate is selected based on the time the algorithm will run; and the iterations have the same complexity than the works in \cite{MJY12,JHA16}. Hence, \cite{NY17,YNW17} generalize and improve the bounds of previous OCO works with static long-term constraints that were motivated by computational complexity of the projections. However, unlike \cite{JHA16}, the learning rate is not adaptive.

The rest of the paper is organized as follows. Sec.\ \ref{sec:model} presents the problem model and Sec.\  \ref{sec:main_results} the main technical results. In Sec.\ \ref{sec:numerical_example}, we present a numerical example where we show the performance of the proposed algorithm depending on $\epsilon$, and compare it to the algorithm in \cite{NY17}. The proof of the main result, Theorem \ref{th:roco}, is in the supplementary material.


\section{Problem Model}
\label{sec:model}

The standard OCO framework can be extended to encompass long-term constraints with additive perturbations as follows. 
Let $C$ be a convex set that contains the admissible or implementable actions, and $g^{(j)}(x) : \R^n \to \R$, $j \in \{1,\dots,m\}$ be  a collection of convex constraints that need to be satisfied on average. Each constraint $g^{(j)}$ has an associated perturbation $b^{(j)}_{t}$ that varies with time. There is no need for the perturbations to be i.i.d., have zero mean, or any other statistical property. The only assumption we will make is that they satisfy a mild Slater condition (see Sec. \ref{sec:assumptions}). In each round $t \in \N$, an agent must select an action $x_t \in C$ without knowledge of the cost function $f_t$ or the perturbation $b_t$. 

We proceed to define the feasible set, the regret, and the constraint violation measure in our problem. To keep notation short, we let $g:= (g^{(1)}, \dots, g^{(m)})$ and $b_t := (b^{(1)}_{t}, \dots, b^{(m)}_{t})$. We define the time-varying feasible set (i.e., the set from which we select the best fixed decision in hindsight)  as
\begin{align}
X_T : = \{ x \in C \mid g(x) +  b_T \preceq 0 \}, \quad \text{where} \quad b^{(j)}_T \in [ \underline{b}^{(j)}_T , \bar b^{(j)}_T], \ j=1,\dots,m
\label{eq:feasible_set}
\end{align}
 with $\bar b^{(j)}_T : = \max \{ b_t^{(j)}, t \in \N \}$, $j \in \{1,\dots,m\}$ and $\underline b_T := \frac{1}{T} \sum_{t=1}^T b_t$.
This is a key difference with previous work where the feasible set is fixed \cite{MJY12} (i.e., $X_T^\text{fixed} : = \{ x \in C \mid g(x) + b \preceq 0 \text{ for some }  b \in \R^m $) or satisfies all the constraints \cite{NY17} (i.e., $X_T^\text{min}: = \{ x \in C \mid g(x) + \bar b_T \preceq 0 \}$). Note that the by construction we have that 
\[
X_T^\text{min} \subseteq X_T \subseteq X_T^\text{max},
\]
where $X_T^\text{max} : = \{ x \in C \mid g(x)  + \underline b_T \preceq 0 \}$ is the set of all the fixed decision that satisfy the constraints on average at time $T$. 
%
The exact value $b_T$ used in the definition of $X_T$ will be specified in Theorem \ref{th:roco}, and will depend on the sequence of perturbations $\{b_t\}$.
%

To avoid confusion with the definition of the regret where the feasible set does not change, we define
\begin{align}
\tilde R(T)  := \sum_{t=1}^T f_t(x_t) - \min_{x \in X_T} \sum_{t=1}^T f_t(x). \label{eq:tilderegret}
\end{align}
%
%
Each action $x_t \in C$ contributes to the aggregate constraint violation
\begin{align}
V(T) & : = \left\| \left[ \sum_{t=1}^{T} g(x_{t}) + b_t \right]^+ \right\|, \label{eq:constraint_violation_def}
\end{align}
where $[z ]^+:= ( \max\{ 0, z^{(1)} \},\dots, \max\{ 0, z^{(m)} \})$ is the projection of each of the components of vector $z \in \R^m$  onto the non-negative orthant. Similar to the regret, we would like that $V(T)$ grows at most sublinearly with $T$ so that $\lim_{T \to \infty}V(T) / T = 0$.  There is no requirement that $V(T) = 0$ for any particular $T \in \N$ or on the rate at which $V(T)$ can grow.
Finally, note that if the sum of the penalties inflicted by a constraint $j \in \{1,\dots,m\}$ is non-positive (i.e., $\sum_{t=1}^T g^{(j)}(x_t) + b^{(j)}_t \preceq 0$), then that constraint does not contribute to the aggregate constraint violation $V(T)$. 
%


\section{Main Results}
\label{sec:main_results}


%
%

\subsection{Proposed algorithm and interpretation}

\begin{algorithm}[t]
\caption{Online Lagrange primal-dual descent/ascent}
\begin{algorithmic}[0]
\STATE {\bfseries Input:} Bregman functions $\psi$ and $\varphi$; vector $f'_1 = 0$; set $C$.
\STATE \textbf{Set:} $x_1 \in C$; $y_1 = 0$; and $\epsilon \in [0,1)$.
\FOR{$t = 1,2,\dots$}
\STATE $\rho \leftarrow 1 / t^\epsilon$
\STATE $ \! \! \! \! \! \! \! \! \! \! \! \ (\circ) \  x_{t+1}  \leftarrow \underset{u \in C}{\arg\min} \{ \mathcal L_{t}(u,y_{t})  +  \frac{1}{\rho} \B (u,x_{t}) \}$ \label{line:x}
\STATE $ \! \! \! \! \! \! \! \! \! \! \! \ (\bullet) \  y_{t+1}  \leftarrow \underset{v \in \R^m_+}{\arg \max } \{ \langle v,  g (x_{t+1}) + b_{t+1} \rangle - \frac{1}{\rho} B_{\varphi} (v,y_{t}) \} $ \label{line:y}
\STATE $f_{t+1} \leftarrow$ play action $x_{t+1}$ and learn cost function
\ENDFOR 
\end{algorithmic}
\label{al:roco}
\end{algorithm}

The main technical contribution of the paper is Algorithm \ref{al:roco}, which allows us to solve the problem presented in Sec.\ \ref{sec:model} with $O(T^\epsilon \vee T^{1-\epsilon})$ regret and $O(T^\epsilon)$ constraint violation.
In short, to handle long-term constraints, we define a Lagrangian-type function
\begin{align}
\mathcal L_{t}(x,y) = \langle f'_{t} (x_{t}) , x \rangle + \langle y, g(x) + b_t \rangle, \label{eq:lagrangian}
\end{align}
where $f'_{t}(x_{t})$ is a (sub)gradient of the cost function in the previous round, and $y \in \R^m_+$ a vector of dual variables. To streamline exposition, in the rest of the paper we will refer to $\mathcal L_{t}(x,y)$ simply as Lagrangian. 
%
Note that the Lagrangian is convex in $x$, concave in $y$, and that it depends on $t$ as the objective function and constraints change in each round.  
The second term of the Lagrangian can be regarded as a penalty or as an adaptive regularizer that allows us to steer the decisions towards the feasible set $X_T$.

Algorithm \ref{al:roco}  is based on a regularized primal-dual proximal gradient method, where we use the general Bregman divergence as the proximal term instead of the usual squared Euclidean distance; see, for example, \cite{BT03}. Recall the Bregman divergence associated with function $\psi$ is defined as
\begin{align}
\B (a,b) = \psi(a) - \psi(b) - \langle a - b, \nabla \psi (b) \rangle, \label{eq:bregman}
\end{align}
where $ \psi$ is usually assumed to be $\sigma_\psi$-strongly convex function and differentiable. 
The primal update $(\circ)$ is equivalent to carrying out an (unconstrained) proximal gradient update with the regularization term $\langle y, g(x) + b_{t} \rangle$.
The regularization or penalty term is updated via the dual update ($\bullet$), which can be regarded as applying a standard proximal gradient ascent since $(g(x_{t+1}) + b_{t+1}) \in \partial_y \mathcal L_{t+1}(x_{t+1},y)$ for a fixed vector $x_{t+1}$.  
%
Interestingly, observe that
\[
\arg \min_{u \in C} \mathcal L_t(u,y_t) = \arg \min_{u \in C} \mathcal \langle f'_{t} (x_{t}) , u \rangle + \langle y_t, g(u)  +b_t \rangle = \arg \min_{u \in C} \mathcal \langle f'_{t} (x_{t}) , u \rangle + \langle y_t, g(u) \rangle ,
\]
and therefore the primal update is oblivious to perturbation $b_t$. Hence, the perturbation is only relevant in our algorithm in the update of the dual variables.\footnote{This is due to the fact that $y$ is the dual variable of the additive perturbation on the constraints. See Sec.\ \ref{sec:sup_ld} in the supplementary material for more details.}

Finally, observe that we use step size $\rho$ equal to $t^{-\epsilon}$ with $\epsilon \in [0,1)$ for both updates, so there is no need to fix in advance the time horizon the algorithm will run. 
Note that when $\epsilon = 0$, then the algorithm corresponds to using constant step size $\rho = 1$.
The algorithm's complexity depends on the structure of the constraints and the functions associated with the Bregman divergence terms in the primal and dual updates. When $g(x)$ is linear (e.g., $g(x) = Ax$) and $\psi,\varphi$ equal to $\frac{1}{2} \|  \cdot \|_2^2$, Algorithm \ref{al:roco} has the same complexity than previous work on OCO with long-term constraints \cite{MJY12,JHA16,YNW17}. In particular, the primal and dual updates can be written as $
x_{t+1}  = \mathcal P_C(x_t - \rho \langle A^T, y_t \rangle)$ and $y_{t+1}  = [y_t + \rho ( g(x_{t+1}) + b_{t+1})]^+$. 

\subsection{Assumptions}
\label{sec:assumptions}
The following are the necessary assumptions to establish the convergence of the proposed algorithm.
\begin{description}
\setlength\itemsep{0em}
\item [Bounded set.] Set $C$ is bounded. There exists a constant $D$ such that $\| u-v \| \le D, \ \forall u,v \in C$. 
\item [Bounded perturbation.] $\| b_t \| < \infty$ for all $t \in \N$. 
\item [Bounded subgradients.] Fix a norm $\| \cdot \|$ and let $\| \cdot\|_*$ denote its dual. There exist constants $F_*$, $G_*$, $G$ such that $\| f'_t(x) \|_* \le F_*$, $\| g(x) + b_t \|_* \le G_*$, $\| g(x) + b_t \| \le G$ for all $x \in C$, $t \in \N$. 
\item [Slater condition.] There exists a  $\eta > 0$ such that $g(\hat x) + b_t + \eta \1 \preceq 0$ for an $\hat x \in C$ and all $t \in \N$. 
\item [Bregman functions.] $\psi$ and $\varphi$ are $\sigma_\psi, \sigma_\varphi $-strongly convex and $L_\psi$, $L_\varphi$-smooth. Also, $\varphi$ is strictly increasing. 
\end{description}

The first assumption is standard in OCO.  
The second and third assumptions are also standard in OCO and ensure that the subgradients of the of primal and dual updates are bounded. 
%
The Slater condition says that there is a set of actions that satisfy the constraints $g(x) + b_t$ strictly for all  $t \in \{1,\dots,T\}$, and is key to ensure that the constraint violation $V(T)$ is sublinear. 
Importantly, the Slater condition assumption is mild in many problems. 
For example, when the perturbation $b_t$ represents the budget available at time $t$, that budget has to be always positive---independently of whether we decide to spend more (i.e., violate the constraint). 
Finally, the assumption that function $\psi$ and $\varphi$ are strongly convex is also standard in the definition of the Bregman divergence. The additional assumption that $\psi$ and $\varphi$ are smooth (hence, $\psi$ and $\varphi$ are upper and lower bounded by a quadratic function) is to streamline exposition in the proofs.\footnote{Technically, all we need is that $\B(u,v)$ is uniformly upper bounded for all $u,v \in C$. Such assumption is also made in previous work and elsewhere to streamline exposition; see, e.g., \cite[Lemma 10]{MJY12}, \cite{DHS11}.} The assumption that $\varphi$ is strictly increasing is necessary to obtain the faster rates on the constraint violation when $\epsilon \in [0,\frac{1}{2})$; note that this is satisfied, for example, by the squared Euclidean distance. 

%



\subsection{Bounds and discussion}
\label{sec:bdisc}


\begin{theorem}
\label{th:roco}
Consider Algorithm~\ref{al:roco} and suppose the assumptions in Sec.\ \ref{sec:assumptions} are satisfied. For any  $X_T : = \{ x \in C \mid g(x) + b_T \preceq 0 \}$ with $b_T  \in  \{ {w  \in [ \underline{b}_T , \bar b_T]} \mid \sum_{t=1}^T \langle y_t , b_{t+1} - w \rangle \le 0\}$,  we have
\begin{align*}
& \tilde R(T) \le \frac{1}{\rho_T}\left(   \frac{L_\psi}{2}D^2 + \frac{L_\varphi}{2}E^2 \right) +  \left( \frac{2F_*^2}{\sigma_\psi}  + \frac{2G_*^2}{\sigma_\varphi} \right)\sum_{t=1}^T \rho_t  = O \left(  T^\epsilon \vee T^{1-\epsilon}  \right)\\
& V(T)  \le G + \frac{L_\varphi }{2\rho_T } E  = O(T^\epsilon)
\end{align*}
where $E$ is a constant that does not depend on $T$ and captures the diameter of the set in which the dual variables are contained. Specifically, 
$
E: = \sqrt {\frac{L_\varphi}{\sigma_\varphi} \left(\frac{2\chi}{\eta}\right)^2 + \frac{2}{\sigma_\varphi}\chi } 
$
where $\chi := \frac{6G^2_*}{\sigma_\varphi}  + 3 F_*D + \frac{L_\psi D^2}{2}$.

\end{theorem}

\underline{\smash{Feasible set}}: The parameter $w  \in [ \underline{b}_T , \bar b_T]$ determines the feasible set $X_T$ used in the definition of the regret in Eq.\ \eqref{eq:tilderegret}. 
Observe that the condition $\sum_{t=1}^T \langle y_t , b_{t+1} - w \rangle \le 0$ always holds for $w = \bar b_T$, since then $b_{t+1}  \preceq \bar b_T := \max\{b_t, t \in \N \}$ and $y_t \succeq 0$ for all $t \in \N$ (this case corresponds to $X_T = X_T^{\text{min}}$ as in \cite{NY17}). If $b_{t+1}$ did not vary over time ($b_t = b$ for all $t \in \N$), then $b_{t+1} = \underline b_T$ and therefore $ \sum_{t=1}^T \langle y_t , b_{t+1} - w \rangle = 0$ (i.e., $X_T = X^{\text{max}}_T$ for all $T \in \N$). Similarly, if $b_t$ were an i.i.d.\ random variable with expected value $b$, then  $\mathbf E(\underline b_T) = b$ and therefore $\mathbf E ( \sum_{t=1}^T \langle y_t , b_{t+1} - w \rangle) \le 0$ for any  $w  \in [ \underline{b}_T , \bar b_T]$.
%
%
It is difficult to characterize $X_T$ for other cases as this does not depend only on the sequences of perturbations $\{b_t\}$ but also on the learning rates (through the dual variables $\{y_t\}$).

\underline{\smash{{Interpretation of the regret bound}}}: The bound has the same structure than the usual OCO bounds.\footnote{See, for example, \cite[Theorem 1]{Zin03}. Or more specifically, the regret bound in the second column on page 4, i.e., $R(T) \le \| F \|^2 \frac{1}{2\eta_T} + \frac{\| \nabla c\|^2}{2} \sum_{t=1}^T \eta_t$.  } For example, for the case where $\psi$ and $\varphi$ are the squared Euclidean distance (i.e., $L_\psi, L_\varphi, \sigma_\psi, \sigma_\varphi = 2$) we have $ \tilde R(T) \le \frac{1}{\rho_T}\left( D^2 + E^2 \right) +  (F_*^2+ G_*^2 )\sum_{t=1}^T \rho_t $. The first term is related to the size of the sets where the primal and dual variables are contained (i.e., $D$ and $E$ respectively) and is inversely proportional to the learning rate at time $T$ (i.e., $\rho^{-1}_{T} = T^\epsilon$). The second term consists of the bounds on the subgradients of the cost functions ($F_*$) and the constraints ($G_*$) multiplied by $\sum_{t=1}^T \rho_t \le 1 + \int_1^T t^{-\epsilon} dt \le 1 + \int_0^T t^{-\epsilon} dt \le 1 + \frac{T^{1-\epsilon}}{1-\epsilon}$.\footnote{We write $O (  T^\epsilon \vee T^{1-\epsilon}  )$ instead of $O (  T^\epsilon \vee \frac{T^{1-\epsilon}}{1-\epsilon}  )$ in Theorem \ref{th:roco} as the interesting range is when $\epsilon \in[0,\frac{1}{2}]$.} 
%
%
The bounds in Theorem \ref{th:roco} are of course useful if the constants are bounded, which is the case for $D$, $F_*$ and $G_*$ by standard OCO assumptions (see Sec.\ \ref{sec:assumptions}). However, for constant $E$ we need more work. To show that this constant exists is the main technical challenge of the paper; we will this discuss it in detail later in the section. 
%


%
Our analysis also allow us to recover the standard OCO bound when the constraints are always satisfied. We have the following corollary to Theorem \ref{th:roco}.

%
%
\begin{corollary}
Suppose that $X_T = C$ (i.e., $g(x) + b_t \preceq 0$ for all $x \in C$ and $t \in \N$). The bound on the regret becomes 
$R(T)  \le \frac{1}{\rho_T} \left( \frac{L_\psi D^2}{2} \right) +  \left( \frac{2 F_*^2}{\sigma_\psi} \right)  \sum_{t=1}^T \rho_t = O(T^\epsilon \vee T^{1-\epsilon})$.
\label{th:ucoro}
\end{corollary}
That is, when the constraints are always satisfied the dual variables will always be equal to zero and therefore $E = 0$ and $G_* = 0$.\footnote{The fact that $G_* = 0$ follows by adding a slack variable $s_t$ to change the inequality constraint to equality, i.e., $g(x) + b_t + s_t = 0$.} 
Hence, by considering perturbed constraints in the learning problem we are adding $ (2 \rho_T)^{-1}L_\varphi E^2 + {2G_*^2}{\sigma_\varphi}^{-1} \sum_{t=1}^T \rho_t $ to the bound of the standard regret in Corollary \ref{th:ucoro}. Such symmetry is not available in previous works \cite{MJY12,JHA16,NY17, YNW17}, and it appears in our work as Algorithm \ref{th:roco} can be regarded, informally, as applying OGD twice (see Lemma \ref{th:bound_regret_first} in the supplementary material for the technical details). As a result, the constants  in the usual OCO bound appear ``duplicated''. 

\underline{\smash{{Interpretation of the constraint violation bound}}}:  The bound on the accumulated constraint violation $V(T)$ consists of two terms. The first term is a constant related to the constraints, and the second term depends on constants $E$ and $L_\varphi$, and are divided by the learning rate at time $T$. Hence, if $\epsilon =  0$ we have that $V(T) \le O(1)$; however, constant constraint violation comes at the price of the regret not being sublinear. Also, observe that for any $\epsilon$ in the range $ [0,\frac{1}{2})$, the constraint violation has better rate than the regret.
\underline{\smash{{Constant $E$}}}: This constant is analogous to constant $D$, which measures the maximum distance between \emph{any} two vectors in the \emph{bounded} set $C$ of primal variables; see Sec.\ \ref{sec:assumptions}. However, we cannot define $E$ in the same way as the dual variables exist in the nonnegative orthant (which is an \emph{unbounded} set).
Instead, we show that the difference between the vectors generated by the dual update in Algorithm \ref{al:roco} is bounded (not \emph{any} two vectors in $\R^m_+$). Or equivalently, that the sequence of dual variables obtained with Algorithm \ref{al:roco} remains bounded for all $t \in \N$; see  Lemma \ref{th:bounded_dual_set} in the supplementary material. 

To ensure that $\| y_t \|$ is bounded for any $t \in \N$, we rely on the Slater condition. In brief, this condition requires that there exists an $x \in C$ such that $g(x) + b_t +\eta \1 \preceq 0$ for some scalar $\eta >0$, and ensures that the dual variables in Algorithm \ref{al:roco} are attracted to a bounded set within $\R^m_+$.\footnote{This type of behavior is typical in \emph{dual} subgradient methods. See, for example, Figure 8.2.6.\ in (Bertsekas et al.\ 2003). This is also discussed in detail in \cite{NO09}; see Lemma 1 in \cite{NO09}.} And since $y_t$ at $t=1$ is bounded, the sequence of dual variables will remain bounded. The technical challenge is to characterize the diameter of the set to which these dual variables are attracted since unlike standard optimization with a fixed objective function, in OCO the cost functions vary over time and, indirectly, the (bounded) sets to which the dual variables are attracted. See Proposition \ref{th:bounded_multipliers} and discussion in Section \ref{sec:bounded_dual_variables} in the supplementary material. 

Finally, we note that  $E = O(D^2)$. Observe that when $\psi$ and $\varphi$ are the squared Euclidean distance, 
$E$ gets simplified to $\sqrt{{2\chi^2 }{\eta^{-2}} + \chi}$ where $\chi:= \frac{6G^2_*}{\sigma_\varphi}  + 3 F_*D + \frac{L_\psi D^2}{2}$. This last observation implies that $\tilde R(T) \le O(D^4 (  T^\epsilon \vee T^{1-\epsilon}  ))$.

\underline{\smash{Constrained convex optimization}}: Our results can also be applied to constrained optimization problems. The following corollary to Theorem \ref{th:roco} establishes the convergence of a constrained convex program with relaxed constraints and primal averaging.

\begin{corollary}
\label{th:coro}
Consider the setup of Theorem \ref{th:roco} where the objective function and constraints are constant (i.e., $f_t = f$ and $b_t =b$ for all $t \in \N$) and step size $\rho_t = \alpha t^{-\epsilon}$ with $\alpha >0$. We have that 
\[
\textup{(i)} \quad f(\bar x_T ) - f^\star  \le O \left(  \frac{1}{\alpha T^{1-\epsilon}} + \frac{\alpha}{T^{\epsilon}}   \right) \qquad  \text{and} \qquad \textup{(ii)} \quad \| [ g ( \bar x_T )  + b ]^+ \| \le O\left( \frac{1}{\alpha T^{1-\epsilon}} \right)
\] 
where $f^\star := \min_{x \in X} f(x)$ with $X= \{x \in C \mid g(x) + b \preceq 0 \}$ and $\bar x_T := \frac{1}{T}\sum_{t=1}^T x_t$. 

\end{corollary}
%
The result recovers the upper bound on the objective and constraint violation in Proposition 1 in \cite{NO09} when $\epsilon = 0$ (fixed step size), but also ensures that $f(\bar x_T) \to f^\star$ and $\bar x_T$ converges to a vector in $X$ asymptotically  as $T \to \infty$ for any $\epsilon \in (0,1)$.


\section{Numerical Example}
\label{sec:numerical_example}

We present a variation of the example in \cite[Sec.\ 5]{YNW17} where the actions made by the decision maker affect the constraints. In short, consider a geo-distributed datacenter that consists of a front-end router and $n$ clusters distributed in different geographical zones. 
Jobs arrive in the front-end router and must be scheduled to one of the clusters. 
The cost of executing jobs depends on the electricity cost of running each of the clusters---which varies across sites as each cluster buys power from its local market. The goal is to schedule jobs to clusters to minimize the total electricity cost while ensuring that all the jobs are served. Importantly, the cost of electricity is not known at the time to schedule jobs.
We  model the problem above as an OCO as follows. Divide time in slots of equal duration and let $x_t \in [0,1]^n$ be the fraction that each of the clusters is utilized at time $t \in \N$.
The cost functions are assumed to be linear (i.e., $f_t (x_t) = \langle l_t, x_t \rangle $ with $l_t \in \R^n_+$) and capture the price of electricity. The constraints are given by $b_t \le \langle a, x_t \rangle $, where $a \in \R^n_+$ captures the efficiency of each cluster (i.e., the number of jobs it can serve per time slot) and $b_t$ the jobs that arrive in the front-end router at time $t$. 

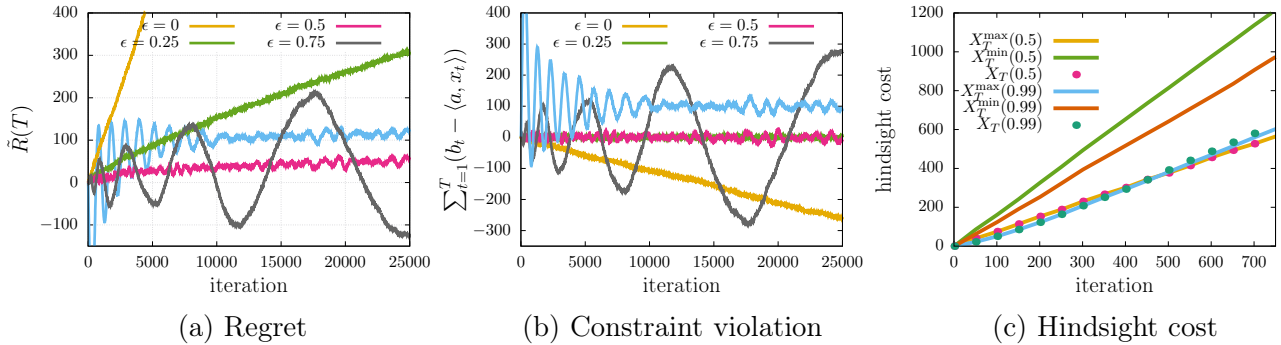
\begin{figure}
\centering
\begin{tabular}{ccc}
\hspace{-1.7em}
{\resizebox{0.35\columnwidth}{!}{\input{Figures/cost.tex}}} & 
\hspace{-1.5em}
{\resizebox{0.35\columnwidth}{!}{\input{Figures/violation.tex}}} &  
\hspace{-1.5em}
{\resizebox{0.35\columnwidth}{!}{\input{Figures/costvs.tex}}} \\
(a) Regret & (b) Constraint violation & (c) Hindsight cost
\end{tabular}
\label{fig:plots}
\caption{{Illustrating the (a) regret and (b) constraint violation of the example in Sec.\ \ref{sec:numerical_example}. The light blue lines in (a) and (b) correspond to the algorithm in \cite{NY17}. Subfigure (c) shows the value of the cost in hindsight depending on the set from which we select the best fixed decision. }}
\end{figure}

We run a simulation with $n =10$ clusters, $\{ l_t \} \in [0,1]^n$, and $b_t = \frac{1}{2} \langle \1 , a \rangle e^{-l_{t-1}(x_{t-1})}$. Hence, the jobs that arrive in a time slot depend on the cost in the previous iteration. 
The simulation results are shown in Fig.\ \ref{fig:plots}, where we evaluate the performance of our algorithm with $\epsilon \in \{0,\frac{1}{4}, \frac{1}{2}, \frac{3}{4}\}$ and compare it with the algorithm in \cite{NY17}  (indicated in blue).  First, observe that by selecting $\epsilon \in \{0,  \frac{1}{4}, \frac{1}{2}\}$, we can trade regret for constraint violation. In particular, with small $\epsilon$, we are sacrificing regret for lower constraint violation---in this example, lower constraint violation means that the jobs have to wait less time at the front-end router before they can be served. When $\epsilon$ is larger than $1/2$, we obtain that our algorithm has a surprising behavior that was not observed in previous work: the regret and constraint violation have a sine-wave form where the period of the wave increases with time.  The latter is shown in the Fig.\ \ref{fig:plots} with $\epsilon = 3/4$, but we observe the same behavior for other values of $\epsilon$ larger than $1/2$. Another interesting observation is the temporal tradeoff between the  regret and constraint violation when $\epsilon = 3/4$; observe from the figure that the peaks of the regret correspond to the lowest values of the constraint violation and vice versa. The lines in light blue show the performance of the algorithm in \cite{NY17},  which compares to our algorithm when $\epsilon = 1/2$ since then both algorithms have $O(\sqrt T)$ regret and constraint violation.\footnote{To compare the results fairly, we select the best policy in hindsight from $X_T$ (the time-varying feasible set in Eq.\ \eqref{eq:feasible_set}) instead of the fixed policy in hindsight that satisfies each constraint individually, i.e., $X_T^{\text{min}}$.} Observe from the figure that our algorithm has significant smaller regret and constraint violation. This result matches the improvements observed in \cite{JHA16} with respect to previous approaches that used a constant learning rate  \cite{MYJ+12}.
Finally, in Fig.\ \ref{fig:plots}c, we show the cost in hindsight depending on whether the best fixed decision is selected from $X_T^{\text{min}}$, $X_T$, or $X_T^{\text{max}}$ with $\epsilon \in \{\frac{1}{2},0.99\}$. Recall that $X_T^{\text{min}} \subseteq X_T \subseteq X_T^{\text{max}}$ and therefore $\min_{x \in X_T^{\text{min}}} \sum_{t=1}^T f_t(x) \ge \min_{x \in X_T^{\text{}}} \sum_{t=1}^T f_t(x) \ge \min_{x \in X_T^{\text{max}}} \sum_{t=1}^T f_t(x). $
Observe from the figure that when $\epsilon = 1/2$, we obtain that the costs with $X_T$ and $X_T^{\text{max}}$ are exactly the same (i.e., the pink dots are exactly on top of the yellow line), whereas when $\epsilon =0.99$, the costs do not coincide exactly ($t \ge 500$). Importantly, notice the large difference between the costs when the best fixed decisions in hindsight are selected from $X_T^{\text{min}}$ and $X_T$. Finally, we note that the sets change over time and are affected by the actions made the decision maker. Observe that when $t= 200$, the cost of the best fixed decision in hindsight is larger with $\epsilon = 0.5$ than with $\epsilon = 0.99$ (i.e., $X_T (0.5) \subset X_T(0.99)$); however, we have the opposite when $t = 700$ (i.e., $X_T (0.99) \subset X_T(0.5)$).
%

%
%


%

\section*{Acknowledgements}
This project has received funding from the European Union’s Horizon 2020 research and innovation programme under the Marie Sk\l{}odowska-Curie grant agreement No 795244.

\bibliography{references}
\bibliographystyle{alpha}

\clearpage

\appendix

\section*{Appendices}

\vspace{1em}

The supplementary material is divided into four sections. Sec.\ \ref{sec:preliminaries} contains the preliminaries and background material. In Sec.\ \ref{sec:technical}, we present the proof of Theorem \ref{th:roco}, and in Sec.\ \ref{sec:proofs} the proofs of the lemmas and propositions in Sec.\ \ref{sec:preliminaries} and \ref{sec:technical}. Finally, in Sec.\ \ref{sec:sup_ld}, we provide additional background material on Lagrange duality and explain how this connects to our approach. 

\section{Preliminaries}
\label{sec:preliminaries}

\subsection{Notation}
We use $\N$, $\R_+$ and $\R^n$ to denote the set of natural numbers, nonnegative real numbers and $n$-dimensional real vectors. The symbol $\1$ indicates the all ones column vector---the dimension of the vector will be implied by the context. We use $f'(x) \in \partial f(x)$ to indicate a subgadient in the subdifferential of $f$ at $x$, and $x_t$ to indicate the $t$'th element in a sequence $\{x_t\}$. For two vectors $u, v \in \R^n$, we write $u \succeq v$ to indicate that $u$ is element-wise larger than or equal to  $v$. The inner product between vectors $u,v$ is indicated with $\langle u,v\rangle$.

\subsection{Bregman divergence}
Let $\B (a,b)$ be the Bregman divergence as defined in Eq.\ \eqref{eq:bregman}. We recall the following well-known identity \cite{CT93}, 
\begin{align}
\B(c,a) + \B (a,b) - \B (c,b)  = \langle \nabla \psi (b) - \nabla \psi (a)  , c-a\rangle, 
\label{eq:tri_property}
\end{align}
which is a generalization of the quadratic identify for the Euclidean norm. It is easy to derive Eq.\ \eqref{eq:tri_property} from the definition of the Bregman divergence in Eq.\ \eqref{eq:bregman}. Observe that
\begin{align*}
& \text{(a)} \quad \B (a,b) = \psi(a) - \psi(b) - \langle a - b, \nabla \psi (b) \rangle, \\
& \text{(b)} \quad \B (c,a) = \psi(c) - \psi(a) - \langle c - a, \nabla \psi (a) \rangle, \\
& \text{(c)} \quad -\B (c,b) = \psi(b) -\psi(c)  - \langle b - c, \nabla \psi (b) \rangle.
\end{align*}

\subsection{Online proximal gradient method }

Let $\phi:\R^n \to \R$ be a convex function. The following update 
\begin{align}
\label{eq:prox_update}
x^+  = \arg \min_{u \in C} \left\{ \phi(u) + \frac{1}{\rho} \B ( u , x )\right\},
\end{align}
corresponds to the standard proximal method where we have replaced the squared Euclidean distance with the Bregman distance.
The following two lemmas are variations of well-known results and correspond to the \emph{proximal} method and \emph{proximal gradient} method (or, mirror-descent); see \cite{Van16}. We state them to measure the progress in one iteration.

\begin{lemma}[One iteration proximal method]
\label{th:one_step_prox}
Consider the proximal update in Eq.\ \eqref{eq:prox_update} where $\phi : \R^n \to \R$ is a convex function, $C \subset \R^n$ a convex set (not necessarily bounded), and $\rho > 0$. For any $z \in C$, the following bound holds
\begin{align}
& \phi (x^+) - \phi (z)   \le \frac{1}{\rho} \left ( \B (z, x)   - \B(z, x^+) -  \B( x^+ ,x)  \right).
\label{eq:one_shot_prox_bound}
\end{align}
\end{lemma}
Obviously, since the bound holds for any $z \in C$, it also holds for the $z$ that minimizes $\phi(z)$. 
Next, we apply the result in Lemma \ref{th:one_step_prox} to the online setting. Let $f$ and $\theta$ be two convex functions from $\R^n \to \R$ and let $\phi (u) = \langle f'(x),  u \rangle + \theta (u)$,
where $ f'(x) \in \partial f (x)$. Function $\theta$ can be regarded as a penalty function or regularizer, and we can ignore it if $\theta(u) = 0$ for all $u \in C$. 
We will relate $\theta$ with the second term of the Lagrangian in Eq.\ \eqref{eq:lagrangian} in the next section.
\begin{lemma}[One iteration proximal gradient method]
\label{th:prox_online_update}
Consider the setup of Lemma \ref{th:one_step_prox} where $\phi (u) = \langle  f'(x),  u \rangle + \theta (u)$ in Eq.\ \eqref{eq:prox_update}.  For any $z \in C$, the following bound holds 
\begin{align*}
f(x) - f(z) + \theta(x^+) - \theta(z)  \le \frac{1}{\rho} \left( \B ( z,x )  - \B ( z,x^+ )\right) +  \frac{2\rho}{\sigma_\psi} \| f'(x) \|^2_*.
\end{align*}
\end{lemma}
The key point from the last lemma is that by only using a subgradient of $f$ we can recover a bound on $f$ itself.
Also, observe that if we let $\theta (x) = 0$ for all $x \in C$; index the objective and the step size with $t$  (i.e., let $f = f_t$ and $\rho = \rho_t$); and sum from $t = 1,\dots, T$, we can follow the rationale of the proof of Theorem 1 in \cite{Zin03}\footnote{Shown also in the proof of Lemma \ref{th:bound_regret_first}.} to recover the standard OCO bound (bound given in Corollary \ref{th:ucoro}). 

\section{Proof of Theorem \ref{th:roco}}
\label{sec:technical}

This section contains the technical results that support the claims in Sec.\  \ref{sec:main_results}. It is divided into two parts. In Sec.\ \ref{sec:regret}, we prove a bound on the regret and constraint violation and show how these depend on the boundedness of the dual variables. Sec.\ \ref{sec:bounded_dual_variables} shows that the dual variables remain uniformly bounded with Algorithm \ref{al:roco}, which is the main technical challenge of the paper. 

\subsection{Regret and constraint violation}
\label{sec:regret}
As explained in Sec.\ \ref{sec:main_results}, the update
\begin{align}
x_{t+1} & = \underset{u \in C}{\arg\min} \left\{ \mathcal L_{t} (u,y_{t}) + \frac{1}{\rho_t} \B (u,x_{t}) \right\} \label{eq:dual_x}
\end{align}
is a generalization of Zinkevich's online gradient descent, and it can be regarded as a proximal gradient update (or, mirror descent) with a regularizer. We can use Lemma \ref{th:prox_online_update} to obtain the following result. 
\begin{lemma}
\label{th:bound_regret_first}
Consider the update in Eq.\ \eqref{eq:dual_x} and let $\{ y_t\}_{t=1}^T$ be an arbitrary sequence of vectors from $\R^m_+$. The following bound holds
\begin{align*}
\tilde R(T)    \le - \sum_{t=1}^T  \langle y_t , g(x_{t+1}) + b_{t+1} \rangle  + \frac{L_\psi}{2\rho_T} D^2 +  \frac{2F_*^2 }{\sigma_\psi}  \sum_{t=1}^T \rho_t
\end{align*}
\end{lemma}
From Lemma \ref{th:bound_regret_first}, we obtain the usual bound on the regret\footnote{See the discussion after Theorem \ref{th:roco}.} with the additional term $- \sum_{t=1}^T \langle y_t , g(x_{t+1}) + b_{t+1} \rangle$ due to the regularizer we have added in the update. To make the term vanish, we can apply a proximal gradient update 
\begin{align}
y_{t+1} & = \underset{v \in \R^m_+}{\arg \max } \left\{ \langle v,  g(x_{t+1})+ b_{t+1} \rangle  - \frac{1}{\rho_t} B_{\varphi} (v,y_{t}) \right\}
\label{eq:dual_y}
\end{align}
since $ \mathcal L_{t} (x_{t+1},y)$ is concave in $y$ for a fixed $x_{t+1}$. The following result is also an application of Lemma \ref{th:prox_online_update}.
\begin{lemma}
\label{th:slackness}
Consider the update in Eq.\ \eqref{eq:dual_y} and suppose there exists a constant $E$ such that $\| y_t \| \le E$ for all $t \in \N$. The following bound holds 
\[
- \sum_{t=1}^T  \langle y_t, g(x_{t+1}) + b_{t+1} \rangle   \le \frac{L_\varphi}{2 \rho_T}  E^2 + \frac{2G_*^2}{\sigma_\varphi} \sum_{t=1}^T \rho_t.
\]
\end{lemma}

Combining Lemmas \ref{th:bound_regret_first} and \ref{th:slackness}, we obtain the bound on the regret in Theorem \ref{th:roco}. It only remains to show that constant $E$ exists and does not depend on $t$. Before we proceed to do that, we establish the bound on the constraint violation. 


\begin{proposition}[Constraint violation] 
Select $x_1 \in C$. The updates in Eq.\ \eqref{eq:dual_x} and \eqref{eq:dual_y} have constraint violation
\label{th:constraint_violation}
\begin{align}
V(T) \le  G +  \frac{L_\varphi}{2\rho_{T}} \| y_{T} \| \label{eq:constraint_violation}
\end{align}
where $\| g(x) + b_t \|\le G$ for all $x \in C$ and $t \in \N$.
\end{proposition}
The proof of the proposition is based as well on OGD arguments. Observe from Eq.\ \eqref{eq:constraint_violation} that if $\| y_t \| \le E$ for all $t \in \N$ (as assumed to obtain the result in Lemma \ref{th:slackness}), then we obtain the claimed bound on the constraint violation in Theorem \ref{th:roco}. We show that  constant $E$ exists in the next section. 

\subsection{Bounded dual variables}
\label{sec:bounded_dual_variables}

We take as starting point a classic result from Lagrange duality in constrained convex optimization, which says that the set of optimal dual variables is \emph{bounded} when the Slater condition holds \cite{Uza58}.
This result is important because when we solve the dual problem with an iterative method, such as the subgradient method,\footnote{See Sec.\ 8.2 in \cite{BNO03} for a detailed explanation of the convergence of the dual subgradient method. See also Lemma 1 and 3 in \cite{NO09}.} we obtain a sequence of dual variables that is attracted to a bounded set. Hence, the dual variables remain bounded for all $t \in \N$.

In our problem, since the Slater condition holds for every constraint (see Sec.\ \ref{sec:assumptions}), we could in principle use the same methodology than in offline constrained convex optimization by defining the time-varying Lagrange dual function
\[
\Psi_t(y) = \min_{u \in C} \left\{ \langle f'_t(x_t), u \rangle + \frac{1}{\rho_t} \B(u,x_t)+ \langle y, g(u) + b_{t+1} \rangle    \right\} .
\]
%
However, that is not possible because each dual function depends on the previous one (through the objective $\langle f'_t(x_t), u \rangle + {\rho_t}^{-1} \B(u,x_t)$), which correlates the set of dual solutions. As a result, we cannot establish that each set $\arg \max_{y \in \R^m_+} \Psi_t(y)$, $t \in \N$ is uniformly bounded. We show this formally in the following proposition. 
%


\begin{proposition}[Bounded optimal dual variables]
\label{th:bounded_multipliers}
 For every $y_t^* \in \arg \max_{y \in \R^m_+} \Psi_t(y)$ and any $\hat x \in C$, we have
\begin{align*}
\| y^*_t \| & \le \frac{1}{\eta} \left(  F_* D + \frac{1}{\rho_t} \left( \B ( \hat x,x_{t} ) - \B (\hat x, x_{t+1}) \right) \right) .
\end{align*}
\end{proposition}
\begin{figure}
\centering
\includegraphics[width=0.65\columnwidth]{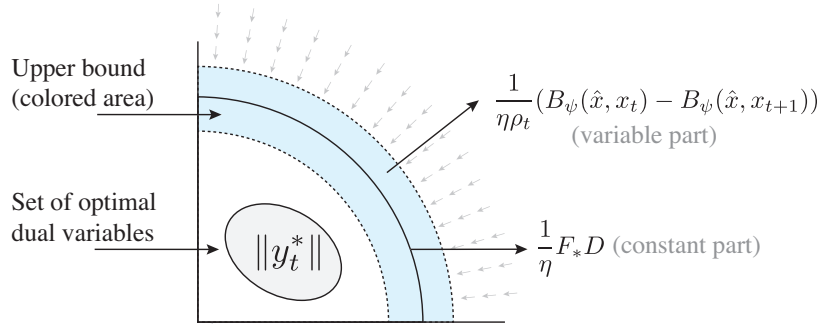}
\caption{Illustrating how the upper bound on the set of dual variables in Proposition \ref{th:bounded_multipliers} depends on the algorithm's parameters. The blue area indicates that the upper bound may vary.  Note that the variable part can be negative.}
\label{eq:balls}
\end{figure}
Proposition \ref{th:bounded_multipliers} gives an upper bound on the set of optimal dual variables for each $t$. The bound consists of two parts: a fixed part $\frac{1}{\eta} F_*D$ and a variable part $\frac{1}{\eta \rho_t} ( \B (  \hat x , x_t ) - \B (\hat x, x_{t+1} ))$ that depends on the previous primal point $x_t$ and step size $\rho_t$. 
The bound is illustrated in  Fig.\ \ref{eq:balls} schematically.
Note that if we upper bound the term $\B (  \hat x , x_t ) - \B (\hat x, x_{t+1} )$ in Proposition \ref{th:bounded_multipliers} by a constant, then the variable term in the upper bound increases with $t$ since $\rho_t \to 0$ as $t\to \infty$ (unless $\epsilon = 0$ and so $\rho_t = 1$). Hence, we cannot claim that each the set of optimal dual variables is uniformly bounded and, therefore, ensure that the sequence of dual variables $\{ y_t\}$  is attracted to a bounded set---which is key to ensure that constant $E$ exists. 
%
%

To deal with the issue mentioned above, we adopt another strategy and show  that the dual variables remain bounded over a span of iterations. Namely, we do not measure the behavior (or, progress) of the dual variables in one step, but over multiple steps.
The intuition behind our strategy is that the variable terms in Proposition \ref{th:bounded_multipliers} ``cancel out'' when we consider the average set of optimal dual variables. 
The length of the interval we consider is proportional to the step size and given by 
\[
S : = \{ t , \dots, t + \lceil t^{\epsilon}  \rceil \} \qquad t  \in \N,
\]
where $\lceil \cdot \rceil$ is the ceiling function.

We proceed to present a lemma that measures the difference of the dual variables according to the Bregman distance over a span of $\lceil t^{
\epsilon} \rceil$ iterations. We present first the following preliminary lemma. 

\begin{lemma}[Bounded subsequences]
\label{th:subs}
Let $\rho_t = t^{-\epsilon}$ and $S : = \{ t, \dots, t + \lceil t^\epsilon \rceil \}$, $t  \in \N$. The following bounds hold
\[
\textup{(i)} \ \ \log(2)\le \sum_{i \in S}  \rho_i \le 3 \qquad \textup{and} \qquad \textup{(ii)} \ \  \sum_{i \in S} \rho_i^2 \le 3.
\]
\end{lemma}

\begin{lemma}
\label{th:attraction}
Consider the setup of Theorem \ref{th:roco}. For any $t \in \N$ and $k \in \{1,\dots,\lceil  t^{\epsilon} \rceil \}$, the following bound holds
\begin{align}
B_\varphi (0,y_{t+k})  - B_\varphi (0, y_{t} )  \le   \chi  - \eta \sum_{i \in S}  \rho_i \| y_t \| 
\label{eq:shift}
\end{align}
where $\chi := \frac{6G^2_*}{\sigma_\varphi}  + 3 F_*D + \frac{L_\psi D^2}{2}$ is a constant that does not depend on $t$.
\end{lemma}
There are two important observations from the last lemma. 

\begin{itemize}
\item \textbf{Observation 1:} When $\eta \sum_{t \in S}  \| y_t \|$ is larger than $\chi$, then the LHS of  Eq.\ \eqref{eq:shift} is negative. Hence,  $y_{t+\lceil t^{\epsilon} \rceil}$ is closer to the origin than $y_t$ with respect to the Bregman divergence ``metric''. Importantly, $\chi$ does not depend on $t$.  

\item \textbf{Observation 2:} For any $t \in \N$, the maximum increment of the dual variables in $\lceil t^{\epsilon} \rceil$ iterations is $\chi$ . That is, $B_\varphi (0,y_{t+\lceil t^{\epsilon} \rceil} )  - B_\varphi (0, y_{t} )  \le   \chi $ for any $t \in \N$ since the second term in the RHS of Eq.\ \eqref{eq:shift} is nonnegative.
\end{itemize}
Using these two observations, we can establish an upper bound on the dual variables for all $t \in \N$. We have the following lemma.

\begin{lemma}[Bounded dual set] 
\label{th:bounded_dual_set}
Consider the setup of Theorem \ref{th:roco} and $\chi$ as defined in Lemma \ref{th:attraction}.
For all $t=1,2,\dots$ we have 
\[
\| y_{t} \|  \le E:= \sqrt {\frac{L_\varphi}{\sigma_\varphi} \left(\frac{2\chi}{\eta}\right)^2 + \frac{2}{\sigma_\varphi}\chi },
\]
\end{lemma}

\section{Proofs}
\label{sec:proofs}

\subsection{Proof of Lemma \ref{th:one_step_prox}}

The update in Eq.\ (\ref{eq:prox_update}) is equivalent to $ x^+ = \arg \min_{u \in \R^n} \{ \phi (u) + \rho^{-1} \B(u,x)  + I_C(u)\}$ where $I_C(u)$ is the indicator function. That is, $I_C(u) = 0 $ if $u \in C$ and $I_C(u) = + \infty$ if $u \notin C$. From the optimality condition,  we have 
\[
0 \in \phi' (x^+)  + \frac{1}{\rho} \nabla \B(x^+,x) +  s,
\] 
where $\phi' (x^+) \in \partial \phi(x^+)$, and $s$ is a vector in the normal cone of $C$ at $x^+$, i.e., $s \in N_C (x^+) := \{ s \in \R^n \mid \langle s , z- x^+  \rangle \le 0, \forall z \in C \}$. Hence, if we multiply the last equation across by $(z - x^+ )$ we obtain
\begin{align}
0 & = \langle \phi' (x^+),  z - x^+  \rangle + \frac{1}{\rho} \langle \nabla \B(x^+, x), z - x^+ \rangle  +  \langle s, z - x^+  \rangle \notag \\
 & \le \langle \phi' (x^+),  z - x^+  \rangle + \frac{1}{\rho} \langle \nabla \B(x^+, x), z - x^+ \rangle. \notag
\end{align}
Next, rearranging terms and using the fact that $\phi(x^+) - \phi(z) \le \langle \phi' (x^+),  x^+ -z  \rangle $ (since $\phi$ is convex) yields
\begin{align*}
\phi(x^+) - \phi(z) 
& \le \frac{1}{\rho} \langle  \nabla \B(x^+,x) , z - x^+ \rangle \\
& =  - \frac{1}{\rho} \langle  -  \nabla \B(x^+,x) ,  z - x^+ \rangle.
\end{align*}
Finally, since $- \nabla \B(x^+,x) = \nabla \psi(x) - \nabla \psi (x^+)$, we can use the identity in Eq.\ (\ref{eq:tri_property}) with $a = x^+$, $b = x$, $c =z$ to obtain $- \langle  \nabla \psi(x) - \nabla \psi (x^+) ,  z - x^+ \rangle =  \B (z,x) -\B(z,x^+) -  \B (x^+,x)$, which concludes the proof.
%


\subsection{Proof of Lemma \ref{th:prox_online_update}}

From Lemma \ref{th:one_step_prox}, we have 
\begin{align}
\langle f' (x) , x^+ - z \rangle + \theta(x^+) - \theta(z)  \le \frac{1}{\rho} \left ( \B (z ,x)   - \B( z , x^+) -  \B( x^+, x)  \right). 
\label{eq:gradprox_simp}
\end{align} 
Add $\langle f'(x),  x - x^+ \rangle$ to both sides of Eq.\ (\ref{eq:gradprox_simp}) and use the fact that $ -  \B( x^+, x) \le - (\sigma_\psi/2) \| x^+ - x\|^2$ (since $\psi$ is strongly convex) to obtain
\begin{align*}
\langle f' (x) , x - z \rangle  + \theta (x^+) - \theta (z)  \le \frac{1}{\rho} \left ( \B ( z ,x)   - \B( z, x^+) \right) + \Delta,
\end{align*} 
where $\Delta :=\langle f'(x),  x - x^+ \rangle - {\sigma_\psi}({2\rho})^{-1} \| x^+ - x\|^2 $. Next, observe that 
\begin{align*}
\Delta
&  \le \sup_{z \in C} \left\{ \langle f'(x),  x - z \rangle - \frac{\sigma_\psi}{2\rho} \|  x - z\|^2 \right\} \\
 &  = \frac{2\rho}{\sigma_\psi}\| f'(x) \|_*^2,
 \end{align*} 
(since the convex conjugate of $\lambda \phi (u)$ for $\lambda > 0$ is $\lambda \phi^*(u^*/\lambda)$; see \cite[pp.\ 475]{RW98}). Finally, since $f(x) - f(z) \le \langle f' (x),  x -z  \rangle $, we obtain the stated result. 

\subsection{Proof of Lemma \ref{th:bound_regret_first}}

Let $\phi(u) = \langle f'_t(x_t),u \rangle + \langle y , g(u) + b_{t+1} \rangle$ in Lemma \ref{th:prox_online_update}. {Note that we use $b_{t+1}$ instead of $b_t$ as the perturbation does not affect the primal variables update (see  the second paragraph in Sec.\ \ref{sec:main_results}). We have
\begin{align}
&  f_{t}(x) -  f_{t}(z) +  \langle y,  g(x^+) + b_{t+1}\rangle  -  \langle y,  g(z) + b_{t+1}\rangle  \\
& \qquad \qquad  \le  \frac{1}{\rho}  \B (z,x) - \B( z,x^+) +   \frac{2F_*^2}{\sigma_\psi} \rho,  \label{eq:teq1}
\end{align}
where we have used the fact that $\| f_t'(x)\|_* \le F_*$ by assumption (see Sec.\ \ref{sec:assumptions}). 
Next, let $x^+ = x_{t+1}$, $x = x_t$, $y = y_t$, $\rho = \rho_t$ in Eq.\ (\ref{eq:teq1}). Summing from $t=1,\dots,T$ and rearranging terms yields
\begin{align*}
\sum_{t=1}^T f_t(x_t) - \sum_{t=1}^T f_t(z)
& \le  \sum_{t=1}^T \frac{1}{\rho_t}  (\B (z, x_t) - \B(z, x_{t+1}) )    +  \frac{2F_*^2}{\sigma_\psi}  \sum_{t=1}^T \rho_t \\
& \quad - \sum_{t=1}^T \langle y_t , g (x_{t+1}) + b_{t+1} \rangle + \sum_{t=1}^T \langle y_t,  g(z) + b_{t+1}\rangle.
\end{align*}
Now, select $z \in X_T$ and observe that
\begin{align}
\sum_{t=1}^T \langle y_t,  g(z) + b_{t+1}\rangle 
& = \sum_{t=1}^T \langle y_t,  g(z) + b_{t+1} + b_T -  b_T \rangle \\
& = \sum_{t=1}^T \langle y_t,  g(z) +  b_T \rangle  + \sum_{t=1}^T \langle y_t,  b_{t+1} - b_T \rangle \\
& \le 0
\end{align}
where the last equation follows since $g(z) + b_T \preceq 0$ for any $z \in X_T \subseteq X^{\text{max}}_{T}$ (by construction); and $ \sum_{t=1}^T \langle y_t,  b_{t+1} - b_T \rangle \le 0$ by the choice of $b_T$ (see Theorem \ref{th:roco}). Hence, 
\[
\tilde R(T)  \le  \sum_{t=1}^T \frac{1}{\rho_t}  (\B (z, x_t) - \B(z, x_{t+1}) )    +  \frac{2F_*^2}{\sigma_\psi}  \sum_{t=1}^T \rho_t  - \sum_{t=1}^T \langle y_t , g (x_{t+1}) + b_{t+1} \rangle .
\]
Finally,  we can upper bound the first term in the RHS of the last equation as follows
\begin{align*}
 \sum_{t=1}^T \frac{1}{\rho_t}  (\B (z, x_t) - \B(z, x_{t+1}) ) 
&  = \frac{\B(z, x_{1})}{\rho_1} -\frac{\B(z, x_{T+1})}{\rho_T}   + \sum_{t=2}^T   \B (z, x_t) \left( \frac{1}{\rho_t} - \frac{1}{\rho_{t-1}}\right) \\
&  \le \frac{\B(z, x_{1})}{\rho_1}   + \sum_{t=2}^T   \B (z, x_t) \left( \frac{1}{\rho_t} - \frac{1}{\rho_{t-1}}\right) \\
&  \stackrel{\text{(a)}}{\le} \frac{L_\psi \| z- x_1\|^2}{2\rho_1}   + \sum_{t=2}^T   \frac{L_\psi \| z- x_t \|^2}{2} \left( \frac{1}{\rho_t} - \frac{1}{\rho_{t-1}}\right) \\
&  \stackrel{\text{(b)}}{\le} \frac{L_\psi D^2}{2\rho_1}   + \sum_{t=2}^T   \frac{L_\psi D^2}{2} \left( \frac{1}{\rho_t} - \frac{1}{\rho_{t-1}}\right) \\
& =    \frac{L_\psi D^2}{2\rho_T},
 \end{align*}
 where (a) follows from the smoothness of $\psi$, and (b) from the assumption that $\max_{u,v \in C} \|u-v\| \le D$ (see Sec.\ \ref{sec:assumptions}).


\subsection{Proof of Lemma \ref{th:slackness}}

Let $\phi(y)  = \langle y , - (g(x_{t+1}) + b_{t+1})\rangle$ in Lemma \ref{th:prox_online_update} and fix $y^+ = y_{t+1}$ and $y = y_t$. Summing from $t=1,\dots,T$ we have 
\begin{align*}
&  \sum_{t=1}^T  \langle y_t, -(g(x_{t+1}) + b_{t+1}) \rangle - \langle z, - (g(x_{t+1}) + b_{t+1})\rangle \\ 
 & \qquad \le \sum_{t=1}^T \frac{1}{\rho_t}  ( B_\varphi (z, y_t)   - B_\varphi (z, y_{t+1}) ) + \frac{2G_*^2}{\sigma_\varphi} \sum_{t=1}^T \rho_t
\end{align*}
where we have used the fact that $\| g(x) + b_t \|_* \le G_*$ for all $x \in C$, $t \in \N$ by assumption; see Sec.\ \ref{sec:assumptions}.  
Next, let $z = 0$ and rearrange terms to obtain
\begin{align*}
 - \sum_{t=1}^T  \langle y_t, g(x_{t+1}) + b_{t+1} \rangle &  \le \frac{1}{\rho_1} {B_\varphi (0, y_1) } - \frac{1}{\rho_T} {B_\varphi (0, y_{T+1})} \\
&  \qquad  + \sum_{t=2}^T B_\varphi (0, y_t)  \left( \frac{1}{\rho_t}  -  \frac{1}{\rho_{t-1}} \right)  + \frac{2G_*^2}{\sigma_\varphi} \sum_{t=1}^T \rho_t.
\end{align*}
Dropping the second term in the RHS of the last equation and using the fact that $B_\varphi (0, y_t) \le 2^{-1} L_\varphi \| y_t \|^2$ (by the smoothness of $\varphi$) and $\| y_t \| \le E$ for all $t \in \N$ (by assumption), we obtain  
\begin{align*}
 - \sum_{t=1}^T  \langle y_t, g(x_{t+1}) + b_{t+1}\rangle  
 & \le  \frac{L_\varphi}{2} E^2  \left( \frac{1}{\rho_1} + \sum_{t=2}^T \left( \frac{1}{\rho_t}  -  \frac{1}{\rho_{t-1}} \right) \right) + \frac{2G_*^2}{\sigma_\varphi} \sum_{t=1}^T \rho_t \\
&  =  \frac{L_\varphi}{2 \rho_T} E^2 + \frac{2G_*^2}{\sigma_\varphi} \sum_{t=1}^T \rho_t,
\end{align*}
which concludes the proof.



\subsection{Proof of Proposition \ref{th:constraint_violation}}

From the optimality condition of  the update in  Eq.\ \eqref{eq:dual_y} we have
\[
0 \in g(x^+) + b_{t+1}  - \frac{1}{\rho} \left( \nabla \varphi (y^+) - \nabla \varphi (y)\right) -s,
\]
where $s$ is a vector in the normal cone of $\R^m_+$ at $y^+$, i.e., $s \in  N_{\R^m_+} (y^+)$.
Rearranging terms and using the fact that $N_{\R^m_+} (y^+) \subseteq \R^m_-$ (i.e., $s \preceq 0$) yields 
\[
g(x^+) + b_{t+1} \preceq \frac{1}{\rho} (\nabla \varphi (y^+)  - \nabla \varphi(y)).
\]
Next, let $x^+ = x_{t+1}$, $y^+ = y_{t+1}$, $y = y_{t}$,  $\rho = \rho_{t}$ and sum the last equation from $t = 1,\dots, T-1$ to obtain
\begin{align*}
 \sum_{t=1}^{T-1} g (x_{t+1}) + b_{t+1}
& \preceq \sum_{t=1}^{T-1} \frac{1}{\rho_{t}} ( \nabla \varphi (y_{t+1}) - \nabla \varphi(y_{t})) \\
& \stackrel{\text{(a)}}{=} \frac{\nabla \varphi(y_{T})}{\rho_{T-1}} - \frac{\nabla \varphi(y_1)}{\rho_1} +  \sum_{t=2}^{T-1} \nabla \varphi(y_t)  \left( \frac{1}{\rho_{t-1}} - \frac{1}{\rho_{t}}  \right)\\
& \stackrel{\text{(b)}}{\preceq} \frac{\nabla \varphi(y_{T})}{\rho_{T-1}}
\end{align*}
where (a) follows by rearranging terms and (b) by dropping the second and third terms in the RHS of (a) since they are nonpositive (note that $\nabla \varphi (y) \succeq 0$ since $\varphi$ is a strictly increasing function and  $\rho_{t+1} \le \rho_{t}$ for all $t\in \N$). Furthermore, observe that we can write
\[
 \left[ \sum_{t=1}^{T-1} g(x_{t+1}) + b_{t+1} \right]^+ \preceq \frac{\nabla \varphi(y_{T})}{\rho_{T-1}}  \preceq \frac{\nabla \varphi(y_{T})}{\rho_{T}}.
\]
Adding $g(x_1) + b_1$ to both sides
\[
 \left[ \sum_{t=1}^{T} g(x_{t}) + b_{t} \right]^+ \preceq g(x_1) + b_1 +  \left[ \sum_{t=1}^{T-1} g(x_{t+1}) + b_{t+1} \right]^+  \preceq g(x_1) + b_1 + \frac{\nabla \varphi(y_{T})}{\rho_{T}},
\]
and therefore
\begin{align*}
\left\| \left[ \sum_{t=1}^{T} g (x_{t}) + b_{t} \right]^+ \right\| \le  \left\|g (x_1) + b_1  + \frac{ \nabla \varphi (y_{T})}{\rho_{T}}   \right\| \le \|g (x_1) + b_1 \| + \frac{1}{\rho_{T}} \|  \nabla \varphi (y_{T}) \| .
\end{align*} 
Finally, if we use the fact that 
$\| \nabla \varphi (y_{T}) \| \le 2^{-1} L_{\varphi} \| y_{T} \|$ since $\varphi$ is $L_\varphi$-smooth (by assumption; see Sec.\ \ref{sec:assumptions}) and $\| g(x) + b_t \| \le G$ for all $t \in \N$, we obtain the stated result. 



\subsection{Proof of Lemma \ref{th:subs}}

We start with claim (i). From the integral test, we have 
$
\delta \le \sum_{t \in S} \rho_t  \le 1 + \delta
$
where 
\begin{align*}
 \delta: = \int_{t}^{t + \lceil t^{\epsilon} \rceil} \frac{1}{x^\epsilon}\ dx = \left[ \frac{x^{1-\epsilon}}{1 - \epsilon} \right]^{t + \lceil t^{\epsilon} \rceil }_t
\end{align*}
Hence, we need to upper and lower bound $\delta$. For the upper bound, observe that  
\[
\left[ \frac{x^{1-\epsilon}}{1 - \epsilon} \right]^{t + \lceil t^{\epsilon} \rceil }_t \le \left[ \frac{x^{1-\epsilon}}{1 - \epsilon} \right]^{t + t^{\epsilon} + 1 }_t = \frac{(t + t^{\epsilon} + 1)^{1-\epsilon} - t^{1-\epsilon}}{1-\epsilon}.
\] 
The equation in the RHS is decreasing in $t$ for a fixed $\epsilon$, but also decreasing in $\epsilon$ for a fixed $t$. Thus, the maximum is attained when $t =1 $ and $\epsilon = 0$. Hence, $\delta \le 2$. 

We continue with the lower bound. Observe that 
\[
\delta  = \left[ \frac{x^{1-\epsilon}}{1 - \epsilon} \right]^{t + \lceil t^{\epsilon} \rceil }_t \ge \left[ \frac{x^{1-\epsilon}}{1 - \epsilon} \right]^{t + t^{\epsilon}  }_t = \frac{(t + t^{\epsilon} )^{1-\epsilon} - t^{1-\epsilon}}{1-\epsilon}.
\] 
For any $t\in \N$, the minimum is attained when $\epsilon = 1$ and equal to 
\[
\lim_{\epsilon \to 1}  \frac{(t + t^{\epsilon})^{1-\epsilon} - t^{1-\epsilon}}{1-\epsilon} = \log(2).
\]
Hence, $\log(2) \le \sum_{t \in S} \rho_t \le 3$ as claimed.

We proceed to show claim (ii). Using again the fact that $\sum_{t \in S} \rho^2_t  \le 1 + \delta$, we can write
\begin{align*}
 \delta: = \int_{t}^{t + \lceil t^{\epsilon} \rceil} \frac{1}{x^{2\epsilon}}\ dx = \left[ \frac{x^{1-2\epsilon}}{1 - 2\epsilon} \right]^{t + \lceil t^{\epsilon} \rceil }_t \le \left[ \frac{x^{1-2\epsilon}}{1 - 2\epsilon} \right]^{t +  t^{\epsilon} +1 }_t \le \frac{(t + t^{\epsilon} + 1)^{1-2\epsilon} - t^{1-2\epsilon}}{1-2\epsilon}.
\end{align*}
Like in the first case, the maximum is attained when $t = 1$ and $\epsilon = 0$ and therefore $\sum_{t \in S} \rho^2_t \le 3$ as claimed.


\subsection{Proof of Lemma \ref{th:attraction}}

From Lemma \ref{th:prox_online_update} with $\phi(y)  = \langle y , -(g(x_{t+1}) + b_{t+1}) \rangle$ and $z = 0$ we have
\begin{align}
\frac{1}{\rho} \left( B_\varphi ( 0, y^+  )  - B_\varphi  ( 0,y)\right) -  \frac{2\rho}{\sigma_\varphi} \| g(x_{t+1}) + b_{t+1} \|^2_* 
\le \langle g (x_{t+1}) + b_{t+1},  y - 0 \rangle . \label{eq:sup_bound_mult_lemma6}
\end{align}
Now observe that
\begin{align}
 &   \langle y ,  g(x_{t+1}) + b_{t+1} \rangle \notag \\
& \ \  =   \langle y, g(x_{t+1}) + b_{t+1} \rangle +  \langle f_t'(x_t),  x_{t+1} -  x_{t+1} \rangle \notag \\
\text{(a)} & \ \ \le   \langle f_t'(x_t),  z \rangle + \langle y, g(z) + b_{t+1} \rangle  - \langle f_t'(x_t),  x_{t+1} \rangle \notag \\ 
& \qquad +  {\rho}^{-1} \left ( \B ( z, x_t)   - \B( z, x_{t+1}) -  \B( x_{t+1} ,x)  \right) \notag \\
& \ \ =   \langle y, g(z) + b_{t+1} \rangle + \langle f_t'(x_t),  z - x_{t+1} \rangle  +  {\rho}^{-1} \left ( \B (z, x_t)   - \B( z, x_{t+1}) -  \B( x_{t+1} ,x_t)  \right) \notag \\
\text{(b)} & \ \  \le   \langle y, g(z) + b_{t+1} \rangle + F_*D +  {\rho}^{-1} \left ( \B ( z, x_t)   - \B( z, x_{t+1})   \right) \notag \\
\text{(c)} & \ \ \le  - \eta \| y \| +  F_*D +  {\rho}^{-1} \left ( \B (  \hat x, x_t)   - \B( \hat x, x_{t+1})   \right) \notag
\end{align}
where (a) follows from Lemma \ref{th:one_step_prox} for any $z \in C$; (b)  since $\langle  f_t'(x_t) , z - x_{t+1} \rangle \le \| f'_t(x_t) \|_* \| z - x_{t+1} \| \le F_*D$ by H\"older's inequality and by dropping $- \rho^{-1} \B( x_{t+1} ,x_t) $; and (c) by letting $z = \hat x$ (a Slater point that satisfies all the constraints; see Sec.\ \ref{sec:assumptions}) and the fact that 
\[
\langle y, g(\hat x) + b_{t+1} \rangle \le -\eta \langle y, \1 \rangle = - \eta \| y \|_1 \le -\eta \| y \|\]
for some $\eta > 0$. 
%
Hence, by multiplying Eq.\ \eqref{eq:sup_bound_mult_lemma6} across by $\rho$ we have
\[
B_\varphi  ( 0, y^+  ) -  B_\varphi  ( 0, y )  
  \le   
     \frac{2\rho^2}{\sigma_\varphi} \| g(x_{t+1}) + b_{t+1}\|^2_*  - \eta \rho \| y \| 
    + \rho F_*D    + \B (\hat x,x_t )  -  \B( \hat x, x_{t+1} ).
 \]
Next, let $y^+ = y_{t+1}$, $y  = y_t$, $\rho = \rho_t$ and sum from $t, \dots, t + k$ with $k \in \N$
\begin{align}
 B_\varphi (0,y_{t + k}  )  - B_\varphi (0,y_{t})    \le   \frac{2G_*^2}{\sigma_\varphi}  \sum_{i = t}^{t+k} \rho_i^2   - \eta  \sum_{i = t}^{t+k} \rho_i \| y_i \| + F_*D  \sum_{i=t}^{t+k} \rho_i  +  \frac{L_\psi D^2}{2}  \label{eq:bdec} 
\end{align}
where in the last equation we have used the fact that $\| g (x) + b_t \|_* \le G_*$ for all $ x \in C$, $t \in \N$  and that $\B (z,x) \le 2^{-1} L_\psi D^2$ for all $x,z \in C$. Finally, by using the bounds in Lemma  \ref{th:subs}, we obtain that 
\begin{align}
 B_\varphi (0,y_{t + k}   )  - B_\varphi (0,y_{t})   
\le   \frac{6G_*^2 }{\sigma_\varphi}   - \eta  \sum_{i = t}^{t+k} \rho_i \| y_i \| + 3 F_*D  +  \frac{L_\psi D^2}{2}  \label{eq:bdec} 
\end{align}
for any $k \in \{1,\dots,\lceil t^\epsilon \rceil\}$.


\subsection{Proof of Lemma \ref{th:bounded_dual_set}}

Define set
\[
\Omega : = \left\{ y \in \R^m_+ \mid \| y \| \le \frac{2\chi}{\eta} \right\},
\]
and consider the following two observations from Lemma \ref{th:attraction}.
\begin{itemize}
\item \textbf{Observation (i).}  For every $y_t \in \R^m_+$ and $k \in \{1,\dots, \lceil t^{\epsilon} \rceil \}$, we have
\[
B_\varphi(0,y_{t + k} ) - B_\varphi(0,y_{t}) \le \chi.
\]
\item \textbf{Observation (ii).}  For every $y_t \in \R^m_+$, if  $y_{t+k} \notin \Omega$ for all $k \in \{1,\dots, \lceil t^{\epsilon} \rceil \}$, then
\begin{align*}
B_\varphi(0,y_{t + \lceil t^{\epsilon} \rceil}) - B_\varphi(0,y_{t}) 
& \le  \chi  - \eta \sum_{i \in S}  \rho_i \| y_i \| \\
& \stackrel{\text{(a)}}{<}    \chi  -  {2\chi}  \sum_{i \in S}  \rho_i \\
& \stackrel{\text{(b)}}{\le} 0
\end{align*}
where (a) follows since  $-\| y \| > - 2 \chi / \eta$ if $y \notin \Omega$, and (b) because $\sum_{i \in S}  \rho_i \ge \log(2) > 1/2$ by Lemma \ref{th:subs}. The last bound also holds if $y_{t + (i-1)\lceil t^\epsilon\rceil + k } \notin \Omega$ for all $i = 1,\dots,M$ with $M \in \{1,2,\dots,\}$, $k \in \{1,\dots, \lceil t^{\epsilon} \rceil \}$. That is, 
\[
 \sum_{i=1}^{M} \left( B_\varphi(0,y_{t + i\lceil t^{\epsilon} \rceil}) - B_\varphi(0,y_{t + (i-1) \lceil t^{\epsilon} \rceil})  \right)  = B_\varphi(0,y_{t + M \lceil t^{\epsilon} \rceil}) - B_\varphi(0,y_{t})
  \le 0.
  \]
  \end{itemize}
Now, consider the case where $y_t \in \Omega$ but $y_{t+1} \notin \Omega$.\footnote{Note that $y_1 \in \Omega$.} Fix $t' = {t + M \lceil t^{\epsilon} \rceil} + k $ for some $M \in \{ 0,1,\dots\}$ and $k \in \{1,\dots, \lceil t^\epsilon \rceil\}$. Combining the two observations, we can write
\[
B_\varphi(0,y_{t'}) \le B_\varphi(0,y_{t}) + \chi.
\]
Next, observe that if we use the fact that $\frac{\sigma_\varphi}{2} \| y \|^2 \le B_\varphi (0,y) \le \frac{L_\varphi}{2} \| y \|^2$ (by the strong convexity and the smoothness of $\varphi$), we can write
\begin{align*}
\| y_{t'} \| & \le \sqrt{ \frac{L_\varphi}{\sigma_\varphi} \| y_t\|^2 + \frac{2}{\sigma_\varphi} \chi} \\
& \stackrel{\text{(a)}}{\le} \sqrt{ \frac{L_\varphi}{\sigma_\varphi} \left( \frac{2\chi}{\eta} \right)^2 + \frac{2}{\sigma_\varphi} \chi} \\
& : = E
\end{align*}
where (a) follows since $y_t \in \Omega$. The last equation concludes the proof since $\max_{y \in \Omega} \| y \| \le E$. 

\subsection{Proof of Proposition \ref{th:bounded_multipliers}}

Let $\phi (z) =   \langle f'_t(x_t), z \rangle +  \langle y^*_t  , g(z) + b_{t+1}  \rangle   $ in Lemma \ref{th:one_step_prox} with $z\in C$ and $y^*_t \in \arg \max_{y \in \R^m_+ } \{ \min_{u \in C} \mathcal L_t (u,y)\}$. We have that
\begin{align*}
&   \langle f'_t(x_t), x^+ -z \rangle +  \langle y^*_t  , g (x^+)  + b_{t+1} - g (z) - b_{t+1} \rangle  \\
& \qquad \le \frac{1}{\rho} \left ( \B(z,x_t)  - \B( z,x^+   ) -  \B ( x^+ , x_t ) \right).
\end{align*}
Rearranging terms and dropping $-  \rho^{-1} \B ( x^+ , x_t )$  yields
\begin{align}
&   - \langle y^*_t, g(z) + b_{t+1} \rangle \notag \\
&  \quad  \le \langle f_t'(x_t), z \rangle  - \langle f_t'(x_t) , x^+ \rangle - \langle y^*_t, g(x^+) + b_{t+1} \rangle   + \frac{1}{\rho} \left ( \B( z,x_t ) - \B( z, x^+  )\right) . \label{eq:basiceqi}
\end{align}
%
Next, observe  that
\begin{align}
 &   \langle f'_t(x_t), z \rangle- \langle f'_t(x_t) , x^+ \rangle - \langle y^*_t, g (x^+)  + b_{t+1} \rangle \notag \\
 & \qquad \le   \langle f'_t(x_t), z \rangle - \min_{u \in C} \{ \langle f'_t(x_t) ,u \rangle + \langle y^*_t, g (u)  + b_{t+1} \rangle \} \notag \\
& \qquad  =  \langle f'_t(x_t), z - x^* \rangle  \label{eq:c_step} \\
&  \qquad  \le F_* D  \label{eq:d_step}
\end{align}
where $x^* \in \arg \min_{u \in C}  \{ \langle f'_t(x_t) ,u \rangle + \langle y^*_t, g (u)  + b_{t+1} \rangle \} $ and Eq.\ \eqref{eq:c_step} follows by complementary slackness \cite[Sec.\ 5.5.2]{BV04}, i.e., 
\[
 \langle f'_t(x_t), x^* \rangle +  \langle y^*_t, g(x^*) + b_{t+1} \rangle 
 =  \langle f'_t(x_t), x^* \rangle
\]
Eq.\ \eqref{eq:d_step} follows since $\langle f'_t(x_t), z - x^* \rangle \le \|  f'_t(x_t) \|_* \| z - x^* \| \le F_* D$ by H\"older's inequality and the fact that set $C$ is bounded (see Sec.\ \ref{sec:assumptions}).
Hence, we can upper bound the first two terms in Eq.\ \eqref{eq:basiceqi} and obtain
\[
- \langle y^*_t, g(z) + b_{t+1}\rangle \le F_* D + \frac{1}{\rho} \left ( \B(z,x )  -  \B( z,x^+ )\right).
\]
Finally, let $z = \hat x$ be a vector that satisfies the Slater condition (see Sec.\ \ref{sec:assumptions}), and note that 
\[
-  \langle y^*_t , g(\hat x) + b_{t+1} \rangle \ge  \eta \langle y^*_t , \1 \rangle =  \eta \| y^*_t \|_1 \ge \eta \| y^* \|.
\]
Using the last bound; diving across by $\eta$; and letting $\rho = \rho_t$, we obtain the stated result.


\section{Lagrange Duality}
\label{sec:sup_ld}

This part is not directly related to the online problem in the main part of the paper, however, we think it may be useful as support material for the readers that are not familiar with Lagrange duality. 
To streamline exposition, we use standard convex optimization notation (e.g., \cite{Roc70,BV04}). 

Let $f_0 : C \to \R$ be a convex function, $f_i: C \to \R$ a collection of $m$ inequality convex constraints\footnote{An equality constraint can be written with two inequality constraints.}, and $C \subseteq \R^n$ a convex set. Define
\[
l(x) := f(x) + \mathbf I (f_1(x),\dots,f_m(x))
\]
where $\mathbf I$ is the indicator function, i.e., $\mathbf I (f_1(x),\dots,f_m(x)) = 0$ if $f_i(x) \le 0$ for all $i=1,\dots, m$, and $\mathbf I (f_1(x),\dots,f_m(x)) = \infty$ if $f_i(x) > 0$ for some $i =1,\dots,m$. Clearly, finding the $x \in C$ that minimizes $l$ is equivalent to finding the $x \in C$ that minimizes $f$ such that every constraint $f_i, i=1,\dots,m$ is satisfied (i.e., the value of every constraint is less than or equal to zero). 

Now, consider the perturbed function
\[
l(x,u) = f(x) + \mathbf I (f_1(x) + u_1,\dots,f_m(x) + u_m)
\]
where $u = (u_1,\dots,u_m)$ is a vector from $\R^m$. Note that $l(x) = l(x,0)$, and that $l(x,u)$ is convex in $u$ for a fixed $x \in C$. Hence, we can write the convex conjugate \cite[Sec.\ 12]{Roc70} of $l(x,u)$ with respect to the perturbation $u$ for a fixed $x \in C$ as follows%
\begin{align*}
l^*(x,u^*) & := \sup_{u \in \R^m} \{ \langle u, u^* \rangle - l(x,u) \} \\
& = \sup_{u \in \R^m} \{ \langle u, u^* \rangle - f(x) - \mathbf I (f_1(x) + u_1,\dots,f_m(x) + u_m) \} \\
 & = - f(x) + \sup_{u \in \R^m} \{ \langle u, u^* \rangle  - \mathbf I (f_1(x) + u_1,\dots,f_m(x) + u_m) \} \\
 & = - f(x) + \sup_{u \in \R^m,  u_i \le - f_i(x),  \ \forall i =1,\dots,m}  \langle u, u^* \rangle 
\end{align*}
where the last equation follows since if $f_i(x) + u_i > 0$ for some $i =1,\dots,m$, then $l^*(x,u^*) = -\infty$ --- which will never be the case since we can always select $u_i  \le -f_i(x)$ as $u_i \in \R$. Next, observe that $\sup_{u \in \R^m,  u_i \le - f_i(x),  \ \forall i =1,\dots,m}  \langle u, u^* \rangle$ can be written as the linear program 
\begin{align*}
\begin{tabular}{llll}
maximize & $\langle u, u^*\rangle$ \\
subject to & $u_i \le - f_i(x)$ & & $ i =1,\dots,m$
\end{tabular}
\end{align*}
It is easy to see that if $u^*_i < 0$ for some $i=1,\dots,m$, then the problem's solution is unbounded above (i.e., we can always select a $u_i$ as negative as we want), and otherwise equal to $- \sum_{i=1}^m u^*_i f_i(x) $. Hence, 
\begin{align*}
l^*(x,u^*)  = 
\begin{cases}
 - f(x) - \sum_{i=1}^m u^*_i f_i(x) & \text{if } u^* \succeq 0  \\
 \infty & \text{otherwise} 
\end{cases}
\end{align*}
Finally, observe that if we let $y = u^*$, we obtain
\[
-l^*(x,y) := \mathcal L(x,y) =  f(x) + \sum_{i=1}^m y_i f_i(x),
\]
which is the classic definition of the Lagrangian for $y \succeq 0$. It is well-known that when the Slater condition holds \cite[Ch.\ 5]{BV04}, then 
\[
\sup_{y \succeq 0} \inf_{x \in C} \mathcal L(x,y) = f^\star = \inf_{x \in C} \sup_{y \succeq 0} \mathcal L(x,y),
\]
where $f^\star = \inf l(x,0)$, i.e., the solution to the ``unperturbed'' problem. Note as well that $\sup_{y \succeq 0} \mathcal L(x,y)$ is indeed like $l(x,0)$; if $x \notin C$ then $\sup_{y\succeq 0} \mathcal L(x,y)$ is equal to $+ \infty$, and otherwise, when $x \in C$, $\sup_{y \succeq 0}\mathcal L(x,y)$ is equal to $f(x)$ .

\textbf{Intuition behind how Lagrange duality fits into our problem.} The perturbed constraints in the main body of the paper can be regarded as if we had the static constraint $f_i(x) = g_i(x) + \bar b_i$ where $\bar b_i = \frac{1}{T} \sum_{t=1}^T b_i^{(t)}$ is the average of the ``perturbations'' for any horizon $T \in \N$. The issue is that the average $\bar b_i$ is not known a priori and is only revealed as we keep playing actions in each round. Hence, we can not use an approach with ``hard'' constraints such as minimizing $l(x,0)$. Instead, we relax the constraints and formulate the Lagrange dual problem. Specifically, we let 
\[
h(y) := \min_{x \in C} \mathcal L(x,y)
\]
and aim to maximize $h$ by carrying out the (sub)gradient ascent update
\[
y^{(t+1)} = [y^{(t)} + \alpha^{(t)} h'(y) ]^+
\]
where $\alpha^{(t)} > 0$ is a step size and $h'(y)$ a subgradient of $h$ at $y$. Note that the projection of the dual variables onto the nonnegative orthant is because $h(y) = -\infty$ if $y_i < 0$ form some $i =1,\dots,m$. The crucial part is that $h'(y)$ is given by $g_i(x^*(y)) + \bar b_i$ where 
\begin{align*}
x^*(y) 
& \in \arg \min_{x \in C} \mathcal L(x,y) \\
& = \arg \min_{x \in C} \left\{ f(x) + \sum_{i=1}^m y_i (g_i(x) + \bar b_i) \right\}\\
& = \arg \min_{x \in C} \left\{ f(x) + \sum_{i=1}^m y_i g_i(x) \right\}
\end{align*}
That is, $x^{*} (y)$ does not depend on the average $\bar b_i$. 

Now observe that we can write the ``noisy'' version of update $y^{(t+1)} = [y^{(t)} + \alpha^{(t)} h'(y) ]^+$ as follows 
\begin{align*}
y^{(t+1)} 
& = [y^{(t)} + \alpha^{(t)}(g(x^*(y)) + \bar b_i + u_i ) ]^+ \\
& = [y^{(t)} + \alpha^{(t)}(g(x^*(y) ) +  b_i) ]^+
\end{align*}
where $u_i = b_i - \bar b_i $ is a ``noise'' vector. In words, the update can be regarded as a ``stochastic'' dual subgradient ascent since $\sum_{t=1}^T u_i^{(t)} = 0$ for all $i = 1,\dots,m$ for any horizon $T$. However, note  that as the average $\bar b_i$ changes, we are changing the set of feasible solutions (or equivalently, the optimization problem itself). This corresponds to the time-varying feasible set $X_T^{\text{max}}$ in the main body of the paper. Recall also that since we do not add any statistical properties to the sequence $\{b_t\}$ of perturbations, we are restricted to comparing our solutions with the more restrictive set $X_T \subseteq X_T^{\text{max}}$. 

The difficulties of applying the approach presented in this section to the online setting are that (i) the cost function varies over time and (ii) that this is not known in advance (i.e., it is learnt after the action has been played). To deal with these issues, we replace $f_t$ for its gradient (as explain in Sec.\ \ref{sec:main_results}) and use a primal-dual proximal gradient approach \cite[Lecture 12]{Van16} as explained in the Sec.\ \ref{sec:preliminaries}.
%

%
%
%

\end{document}

%% file: Figures/cost.tex
\begingroup
  \makeatletter
  \providecommand\color[2][]{%
    \GenericError{(gnuplot) \space\space\space\@spaces}{%
      Package color not loaded in conjunction with
      terminal option `colourtext'%
    }{See the gnuplot documentation for explanation.%
    }{Either use 'blacktext' in gnuplot or load the package
      color.sty in LaTeX.}%
    \renewcommand\color[2][]{}%
  }%
  \providecommand\includegraphics[2][]{%
    \GenericError{(gnuplot) \space\space\space\@spaces}{%
      Package graphicx or graphics not loaded%
    }{See the gnuplot documentation for explanation.%
    }{The gnuplot epslatex terminal needs graphicx.sty or graphics.sty.}%
    \renewcommand\includegraphics[2][]{}%
  }%
  \providecommand\rotatebox[2]{#2}%
  \@ifundefined{ifGPcolor}{%
    \newif\ifGPcolor
    \GPcolorfalse
  }{}%
  \@ifundefined{ifGPblacktext}{%
    \newif\ifGPblacktext
    \GPblacktexttrue
  }{}%
  \let\gplgaddtomacro\g@addto@macro
  \gdef\gplbacktext{}%
  \gdef\gplfronttext{}%
  \makeatother
  \ifGPblacktext
    \def\colorrgb#1{}%
    \def\colorgray#1{}%
  \else
    \ifGPcolor
      \def\colorrgb#1{\color[rgb]{#1}}%
      \def\colorgray#1{\color[gray]{#1}}%
      \expandafter\def\csname LTw\endcsname{\color{white}}%
      \expandafter\def\csname LTb\endcsname{\color{black}}%
      \expandafter\def\csname LTa\endcsname{\color{black}}%
      \expandafter\def\csname LT0\endcsname{\color[rgb]{1,0,0}}%
      \expandafter\def\csname LT1\endcsname{\color[rgb]{0,1,0}}%
      \expandafter\def\csname LT2\endcsname{\color[rgb]{0,0,1}}%
      \expandafter\def\csname LT3\endcsname{\color[rgb]{1,0,1}}%
      \expandafter\def\csname LT4\endcsname{\color[rgb]{0,1,1}}%
      \expandafter\def\csname LT5\endcsname{\color[rgb]{1,1,0}}%
      \expandafter\def\csname LT6\endcsname{\color[rgb]{0,0,0}}%
      \expandafter\def\csname LT7\endcsname{\color[rgb]{1,0.3,0}}%
      \expandafter\def\csname LT8\endcsname{\color[rgb]{0.5,0.5,0.5}}%
    \else
      \def\colorrgb#1{\color{black}}%
      \def\colorgray#1{\color[gray]{#1}}%
      \expandafter\def\csname LTw\endcsname{\color{white}}%
      \expandafter\def\csname LTb\endcsname{\color{black}}%
      \expandafter\def\csname LTa\endcsname{\color{black}}%
      \expandafter\def\csname LT0\endcsname{\color{black}}%
      \expandafter\def\csname LT1\endcsname{\color{black}}%
      \expandafter\def\csname LT2\endcsname{\color{black}}%
      \expandafter\def\csname LT3\endcsname{\color{black}}%
      \expandafter\def\csname LT4\endcsname{\color{black}}%
      \expandafter\def\csname LT5\endcsname{\color{black}}%
      \expandafter\def\csname LT6\endcsname{\color{black}}%
      \expandafter\def\csname LT7\endcsname{\color{black}}%
      \expandafter\def\csname LT8\endcsname{\color{black}}%
    \fi
  \fi
    \setlength{\unitlength}{0.0500bp}%
    \ifx\gptboxheight\undefined%
      \newlength{\gptboxheight}%
      \newlength{\gptboxwidth}%
      \newsavebox{\gptboxtext}%
    \fi%
    \setlength{\fboxrule}{0.5pt}%
    \setlength{\fboxsep}{1pt}%
\begin{picture}(5668.00,3968.00)%
    \gplgaddtomacro\gplbacktext{%
      \csname LTb\endcsname
      \put(946,981){\makebox(0,0)[r]{\strut{}$-100$}}%
      \csname LTb\endcsname
      \put(946,1534){\makebox(0,0)[r]{\strut{}$0$}}%
      \csname LTb\endcsname
      \put(946,2087){\makebox(0,0)[r]{\strut{}$100$}}%
      \csname LTb\endcsname
      \put(946,2640){\makebox(0,0)[r]{\strut{}$200$}}%
      \csname LTb\endcsname
      \put(946,3194){\makebox(0,0)[r]{\strut{}$300$}}%
      \csname LTb\endcsname
      \put(946,3747){\makebox(0,0)[r]{\strut{}$400$}}%
      \csname LTb\endcsname
      \put(1078,484){\makebox(0,0){\strut{}$0$}}%
      \csname LTb\endcsname
      \put(1917,484){\makebox(0,0){\strut{}$5000$}}%
      \csname LTb\endcsname
      \put(2755,484){\makebox(0,0){\strut{}$10000$}}%
      \csname LTb\endcsname
      \put(3594,484){\makebox(0,0){\strut{}$15000$}}%
      \csname LTb\endcsname
      \put(4432,484){\makebox(0,0){\strut{}$20000$}}%
      \csname LTb\endcsname
      \put(5271,484){\makebox(0,0){\strut{}$25000$}}%
    }%
    \gplgaddtomacro\gplfronttext{%
      \csname LTb\endcsname
      \put(198,2225){\rotatebox{-270}{\makebox(0,0){\strut{}\Large $\tilde R(T)$}}}%
      \put(3174,154){\makebox(0,0){\strut{}\Large iteration}}%
      \csname LTb\endcsname
      \put(2266,3574){\makebox(0,0)[r]{\strut{}$\epsilon = 0$}}%
      \csname LTb\endcsname
      \put(2266,3354){\makebox(0,0)[r]{\strut{}$\epsilon = 0.25$}}%
      \csname LTb\endcsname
      \put(4177,3574){\makebox(0,0)[r]{\strut{}$\epsilon = 0.5$}}%
      \csname LTb\endcsname
      \put(4177,3354){\makebox(0,0)[r]{\strut{}$\epsilon = 0.75$}}%
    }%
    \gplbacktext
    \put(0,0){\includegraphics{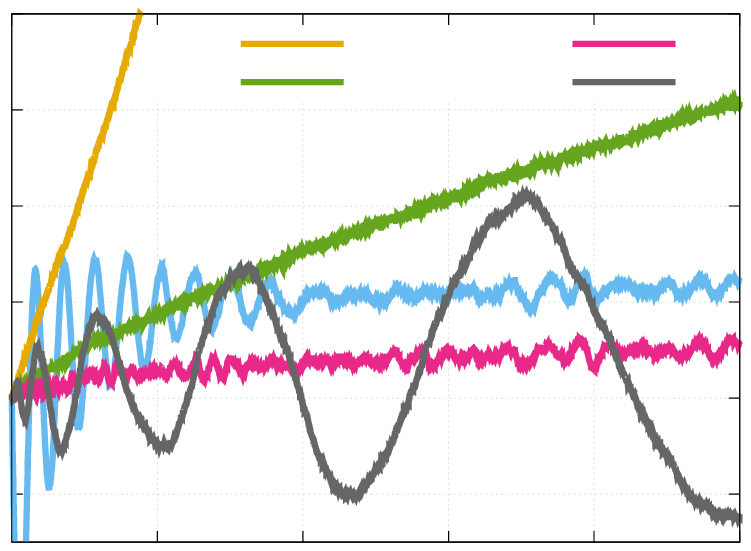}}%
    \gplfronttext
  \end{picture}%
\endgroup

%% file: Figures/violation.tex
\begingroup
  \makeatletter
  \providecommand\color[2][]{%
    \GenericError{(gnuplot) \space\space\space\@spaces}{%
      Package color not loaded in conjunction with
      terminal option `colourtext'%
    }{See the gnuplot documentation for explanation.%
    }{Either use 'blacktext' in gnuplot or load the package
      color.sty in LaTeX.}%
    \renewcommand\color[2][]{}%
  }%
  \providecommand\includegraphics[2][]{%
    \GenericError{(gnuplot) \space\space\space\@spaces}{%
      Package graphicx or graphics not loaded%
    }{See the gnuplot documentation for explanation.%
    }{The gnuplot epslatex terminal needs graphicx.sty or graphics.sty.}%
    \renewcommand\includegraphics[2][]{}%
  }%
  \providecommand\rotatebox[2]{#2}%
  \@ifundefined{ifGPcolor}{%
    \newif\ifGPcolor
    \GPcolorfalse
  }{}%
  \@ifundefined{ifGPblacktext}{%
    \newif\ifGPblacktext
    \GPblacktexttrue
  }{}%
  \let\gplgaddtomacro\g@addto@macro
  \gdef\gplbacktext{}%
  \gdef\gplfronttext{}%
  \makeatother
  \ifGPblacktext
    \def\colorrgb#1{}%
    \def\colorgray#1{}%
  \else
    \ifGPcolor
      \def\colorrgb#1{\color[rgb]{#1}}%
      \def\colorgray#1{\color[gray]{#1}}%
      \expandafter\def\csname LTw\endcsname{\color{white}}%
      \expandafter\def\csname LTb\endcsname{\color{black}}%
      \expandafter\def\csname LTa\endcsname{\color{black}}%
      \expandafter\def\csname LT0\endcsname{\color[rgb]{1,0,0}}%
      \expandafter\def\csname LT1\endcsname{\color[rgb]{0,1,0}}%
      \expandafter\def\csname LT2\endcsname{\color[rgb]{0,0,1}}%
      \expandafter\def\csname LT3\endcsname{\color[rgb]{1,0,1}}%
      \expandafter\def\csname LT4\endcsname{\color[rgb]{0,1,1}}%
      \expandafter\def\csname LT5\endcsname{\color[rgb]{1,1,0}}%
      \expandafter\def\csname LT6\endcsname{\color[rgb]{0,0,0}}%
      \expandafter\def\csname LT7\endcsname{\color[rgb]{1,0.3,0}}%
      \expandafter\def\csname LT8\endcsname{\color[rgb]{0.5,0.5,0.5}}%
    \else
      \def\colorrgb#1{\color{black}}%
      \def\colorgray#1{\color[gray]{#1}}%
      \expandafter\def\csname LTw\endcsname{\color{white}}%
      \expandafter\def\csname LTb\endcsname{\color{black}}%
      \expandafter\def\csname LTa\endcsname{\color{black}}%
      \expandafter\def\csname LT0\endcsname{\color{black}}%
      \expandafter\def\csname LT1\endcsname{\color{black}}%
      \expandafter\def\csname LT2\endcsname{\color{black}}%
      \expandafter\def\csname LT3\endcsname{\color{black}}%
      \expandafter\def\csname LT4\endcsname{\color{black}}%
      \expandafter\def\csname LT5\endcsname{\color{black}}%
      \expandafter\def\csname LT6\endcsname{\color{black}}%
      \expandafter\def\csname LT7\endcsname{\color{black}}%
      \expandafter\def\csname LT8\endcsname{\color{black}}%
    \fi
  \fi
    \setlength{\unitlength}{0.0500bp}%
    \ifx\gptboxheight\undefined%
      \newlength{\gptboxheight}%
      \newlength{\gptboxwidth}%
      \newsavebox{\gptboxtext}%
    \fi%
    \setlength{\fboxrule}{0.5pt}%
    \setlength{\fboxsep}{1pt}%
\begin{picture}(5668.00,3968.00)%
    \gplgaddtomacro\gplbacktext{%
      \csname LTb\endcsname
      \put(946,907){\makebox(0,0)[r]{\strut{}$-300$}}%
      \put(946,1313){\makebox(0,0)[r]{\strut{}$-200$}}%
      \put(946,1718){\makebox(0,0)[r]{\strut{}$-100$}}%
      \put(946,2124){\makebox(0,0)[r]{\strut{}$0$}}%
      \put(946,2530){\makebox(0,0)[r]{\strut{}$100$}}%
      \put(946,2936){\makebox(0,0)[r]{\strut{}$200$}}%
      \put(946,3341){\makebox(0,0)[r]{\strut{}$300$}}%
      \put(946,3747){\makebox(0,0)[r]{\strut{}$400$}}%
      \put(1078,484){\makebox(0,0){\strut{}$0$}}%
      \put(1917,484){\makebox(0,0){\strut{}$5000$}}%
      \put(2755,484){\makebox(0,0){\strut{}$10000$}}%
      \put(3594,484){\makebox(0,0){\strut{}$15000$}}%
      \put(4432,484){\makebox(0,0){\strut{}$20000$}}%
      \put(5271,484){\makebox(0,0){\strut{}$25000$}}%
    }%
    \gplgaddtomacro\gplfronttext{%
      \csname LTb\endcsname
      \put(198,2225){\rotatebox{-270}{\makebox(0,0){\strut{}\Large $\sum_{t=1}^T (b_t - \langle a, x_t \rangle)$}}}%
      \put(3174,154){\makebox(0,0){\strut{}\Large iteration}}%
      \csname LTb\endcsname
      \put(2266,3574){\makebox(0,0)[r]{\strut{}$\epsilon = 0$}}%
      \csname LTb\endcsname
      \put(2266,3354){\makebox(0,0)[r]{\strut{}$\epsilon = 0.25$}}%
      \csname LTb\endcsname
      \put(4177,3574){\makebox(0,0)[r]{\strut{}$\epsilon = 0.5$}}%
      \csname LTb\endcsname
      \put(4177,3354){\makebox(0,0)[r]{\strut{}$\epsilon = 0.75$}}%
    }%
    \gplbacktext
    \put(0,0){\includegraphics{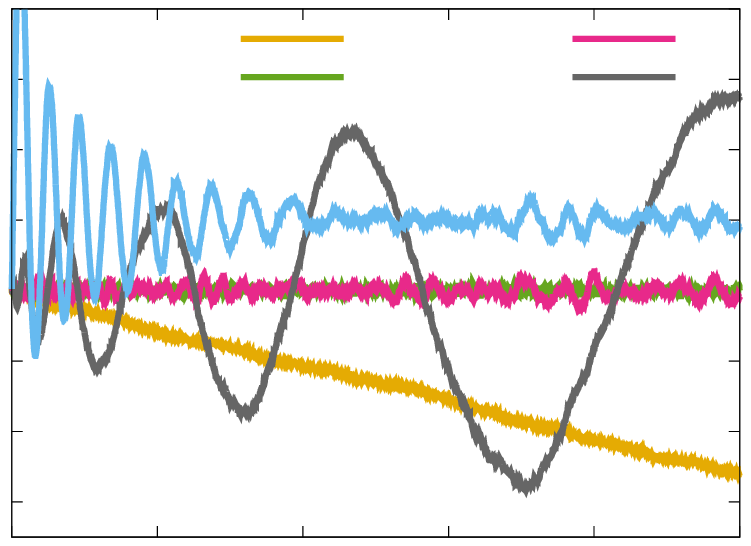}}%
    \gplfronttext
  \end{picture}%
\endgroup

%% file: Figures/costvs.tex
\begingroup
  \makeatletter
  \providecommand\color[2][]{%
    \GenericError{(gnuplot) \space\space\space\@spaces}{%
      Package color not loaded in conjunction with
      terminal option `colourtext'%
    }{See the gnuplot documentation for explanation.%
    }{Either use 'blacktext' in gnuplot or load the package
      color.sty in LaTeX.}%
    \renewcommand\color[2][]{}%
  }%
  \providecommand\includegraphics[2][]{%
    \GenericError{(gnuplot) \space\space\space\@spaces}{%
      Package graphicx or graphics not loaded%
    }{See the gnuplot documentation for explanation.%
    }{The gnuplot epslatex terminal needs graphicx.sty or graphics.sty.}%
    \renewcommand\includegraphics[2][]{}%
  }%
  \providecommand\rotatebox[2]{#2}%
  \@ifundefined{ifGPcolor}{%
    \newif\ifGPcolor
    \GPcolorfalse
  }{}%
  \@ifundefined{ifGPblacktext}{%
    \newif\ifGPblacktext
    \GPblacktexttrue
  }{}%
  \let\gplgaddtomacro\g@addto@macro
  \gdef\gplbacktext{}%
  \gdef\gplfronttext{}%
  \makeatother
  \ifGPblacktext
    \def\colorrgb#1{}%
    \def\colorgray#1{}%
  \else
    \ifGPcolor
      \def\colorrgb#1{\color[rgb]{#1}}%
      \def\colorgray#1{\color[gray]{#1}}%
      \expandafter\def\csname LTw\endcsname{\color{white}}%
      \expandafter\def\csname LTb\endcsname{\color{black}}%
      \expandafter\def\csname LTa\endcsname{\color{black}}%
      \expandafter\def\csname LT0\endcsname{\color[rgb]{1,0,0}}%
      \expandafter\def\csname LT1\endcsname{\color[rgb]{0,1,0}}%
      \expandafter\def\csname LT2\endcsname{\color[rgb]{0,0,1}}%
      \expandafter\def\csname LT3\endcsname{\color[rgb]{1,0,1}}%
      \expandafter\def\csname LT4\endcsname{\color[rgb]{0,1,1}}%
      \expandafter\def\csname LT5\endcsname{\color[rgb]{1,1,0}}%
      \expandafter\def\csname LT6\endcsname{\color[rgb]{0,0,0}}%
      \expandafter\def\csname LT7\endcsname{\color[rgb]{1,0.3,0}}%
      \expandafter\def\csname LT8\endcsname{\color[rgb]{0.5,0.5,0.5}}%
    \else
      \def\colorrgb#1{\color{black}}%
      \def\colorgray#1{\color[gray]{#1}}%
      \expandafter\def\csname LTw\endcsname{\color{white}}%
      \expandafter\def\csname LTb\endcsname{\color{black}}%
      \expandafter\def\csname LTa\endcsname{\color{black}}%
      \expandafter\def\csname LT0\endcsname{\color{black}}%
      \expandafter\def\csname LT1\endcsname{\color{black}}%
      \expandafter\def\csname LT2\endcsname{\color{black}}%
      \expandafter\def\csname LT3\endcsname{\color{black}}%
      \expandafter\def\csname LT4\endcsname{\color{black}}%
      \expandafter\def\csname LT5\endcsname{\color{black}}%
      \expandafter\def\csname LT6\endcsname{\color{black}}%
      \expandafter\def\csname LT7\endcsname{\color{black}}%
      \expandafter\def\csname LT8\endcsname{\color{black}}%
    \fi
  \fi
    \setlength{\unitlength}{0.0500bp}%
    \ifx\gptboxheight\undefined%
      \newlength{\gptboxheight}%
      \newlength{\gptboxwidth}%
      \newsavebox{\gptboxtext}%
    \fi%
    \setlength{\fboxrule}{0.5pt}%
    \setlength{\fboxsep}{1pt}%
\begin{picture}(5668.00,3968.00)%
    \gplgaddtomacro\gplbacktext{%
      \csname LTb\endcsname
      \put(946,704){\makebox(0,0)[r]{\strut{}$0$}}%
      \put(946,1211){\makebox(0,0)[r]{\strut{}$200$}}%
      \put(946,1718){\makebox(0,0)[r]{\strut{}$400$}}%
      \put(946,2226){\makebox(0,0)[r]{\strut{}$600$}}%
      \put(946,2733){\makebox(0,0)[r]{\strut{}$800$}}%
      \put(946,3240){\makebox(0,0)[r]{\strut{}$1000$}}%
      \put(946,3747){\makebox(0,0)[r]{\strut{}$1200$}}%
      \put(1078,484){\makebox(0,0){\strut{}$0$}}%
      \put(1637,484){\makebox(0,0){\strut{}$100$}}%
      \put(2196,484){\makebox(0,0){\strut{}$200$}}%
      \put(2755,484){\makebox(0,0){\strut{}$300$}}%
      \put(3314,484){\makebox(0,0){\strut{}$400$}}%
      \put(3873,484){\makebox(0,0){\strut{}$500$}}%
      \put(4432,484){\makebox(0,0){\strut{}$600$}}%
      \put(4991,484){\makebox(0,0){\strut{}$700$}}%
    }%
    \gplgaddtomacro\gplfronttext{%
      \csname LTb\endcsname
      \put(198,2225){\rotatebox{-270}{\makebox(0,0){\strut{}\Large hindsight cost}}}%
      \put(3174,154){\makebox(0,0){\strut{}\Large iteration}}%
      \csname LTb\endcsname
      \put(2243,3383){\makebox(0,0)[r]{\strut{}$X_T^{\text{max}}(0.5)$}}%
      \csname LTb\endcsname
      \put(2243,3163){\makebox(0,0)[r]{\strut{}$X_T^{\text{min}}(0.5)$}}%
      \csname LTb\endcsname
      \put(2243,2943){\makebox(0,0)[r]{\strut{}$X_T(0.5)$}}%
      \csname LTb\endcsname
      \put(2243,2723){\makebox(0,0)[r]{\strut{}$X_T^{\text{max}}(0.99)$}}%
      \csname LTb\endcsname
      \put(2243,2503){\makebox(0,0)[r]{\strut{}$X_T^{\text{min}}(0.99)$}}%
      \csname LTb\endcsname
      \put(2243,2283){\makebox(0,0)[r]{\strut{}$X_T(0.99)$}}%
    }%
    \gplbacktext
    \put(0,0){\includegraphics{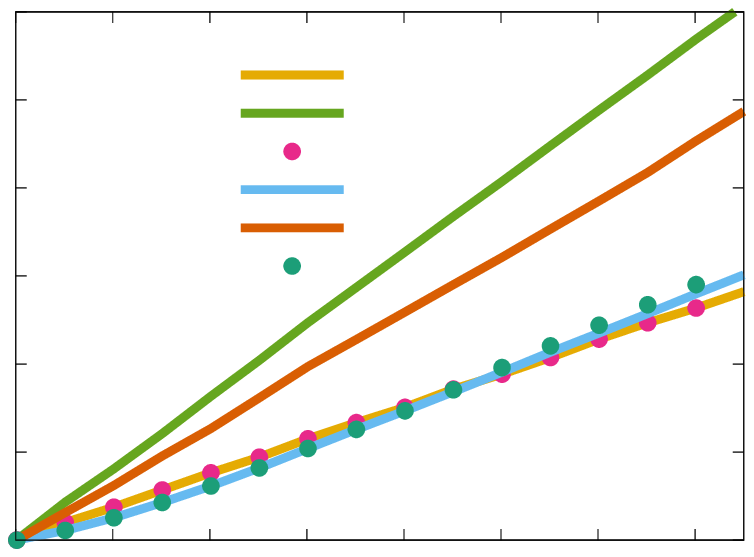}}%
    \gplfronttext
  \end{picture}%
\endgroup

%% file: oco_perturbed.bbl
\newcommand{\etalchar}[1]{$^{#1}$}
\begin{thebibliography}{MTYY09}

\bibitem[AD15]{AD15}
Shipra Agrawal and Nikhil~R. Devanur.
\newblock Fast algorithms for online stochastic convex programming.
\newblock In {\em Proceedings of the Twenty-sixth Annual ACM-SIAM Symposium on
  Discrete Algorithms}, pages 1405--1424, 2015.

\bibitem[BNO03]{BNO03}
Dimitri~P. Bertsekas, Angelia Nedi\'c, and Asuman~E. Ozdaglar.
\newblock {\em Convex analysis and optimization}.
\newblock Athena Scientific, 2003.

\bibitem[BT03]{BT03}
Amir Beck and Marc Teboulle.
\newblock Mirror descent and nonlinear projected subgradient methods for convex
  optimization.
\newblock {\em Operations Research Letters}, 31(3):167 -- 175, 2003.

\bibitem[BV04]{BV04}
Stephne Boyd and Lieven Vandenberghe.
\newblock {\em Convex Optimization}.
\newblock Cambridge University Press, 2004.

\bibitem[CBL06]{CL06}
Nicolo Cesa-Bianchi and Gabor Lugosi.
\newblock {\em Prediction, Learning, and Games}.
\newblock Cambridge University Press, 2006.

\bibitem[CGP16]{CGP16}
Andrew Cotter, Maya Gupta, and Jan Pfeifer.
\newblock A light touch for heavily constrained sgd.
\newblock In {\em 29th Annual Conference on Learning Theory}, volume~49 of {\em
  Proceedings of Machine Learning Research}, pages 729--771, 2016.

\bibitem[CT93]{CT93}
Gong Chen and Marc Teboulle.
\newblock Convergence analysis of a proximal-like minimization algorithm using
  bregman functions.
\newblock {\em SIAM Journal on Optimization}, 3(3):538--543, 1993.

\bibitem[DHS11]{DHS11}
John Duchi, Elad Hazan, and Yoram Singer.
\newblock Adaptive subgradient methods for online learning and stochastic
  optimization.
\newblock {\em Journal of Machine Learning Research}, 12(Jul):2121--2159, 2011.

\bibitem[GIN10]{GIN19}
Eugene Gorbatov, Canturk Isci, and Ripal Nathuji.
\newblock Power efficient resource allocation in data centers, October~19 2010.
\newblock US Patent 7,818,594.

\bibitem[GNT06]{GNT06}
Leonidas Georgiadis, Michael~J. Neely, and Leandros Tassiulas.
\newblock {\em Resource Allocation and Cross-Layer Control in Wireless
  Networks}.
\newblock Foundatations and Trends in Optimization, 2006.

\bibitem[Haz16]{Haz16}
Elad Hazan.
\newblock {\em Introduction to online convex optimization}.
\newblock Now Publishers, Inc., 2016.

\bibitem[JHA16]{JHA16}
Rodolphe Jenatton, Jim~C. Huang, and Cedric Archambeau.
\newblock Adaptive algorithms for online convex optimization with long-term
  constraints.
\newblock In {\em International Conference on Machine Learning}, 2016.

\bibitem[Mey08]{Mey08}
Sean Meyn.
\newblock {\em Control techniques for complex networks}.
\newblock Cambridge University Press, 2008.

\bibitem[MJY12]{MJY12}
Mehrdad Mahdavi, Rong Jin, and Tianbao Yang.
\newblock Trading regret for efficiency: Online convex optimization with long
  term constraints.
\newblock {\em Journal of Machine Learning Research}, pages 2503--2528, Sep
  2012.

\bibitem[MTYY09]{MTY09}
Shie Mannor, John~N. Tsitsiklis, and Jia Yuan~Yu.
\newblock Online learning with sample path constraints.
\newblock {\em Journal of Machine Learning Research}, pages 569--590, Mar 2009.

\bibitem[MYJ{\etalchar{+}}12]{MYJ+12}
Mehrdad Mahdavi, Tianbao Yang, Rong Jin, Shenghuo Zhu, and Jinfeng Yi.
\newblock Stochastic gradient descent with only one projection.
\newblock In {\em Advances in Neural Information Processing Systems 25}, pages
  494--502, 2012.

\bibitem[Nem94]{Nem94}
Arkadi Nemirovski.
\newblock Information-based complexity of convex programming.
\newblock Technical report, Technion, 1994.

\bibitem[NO09]{NO09}
Angelia Nedi{\'c} and Asuman Ozdaglar.
\newblock Approximate primal solutions and rate analysis for dual subgradient
  methods.
\newblock {\em SIAM Journal on Optimization}, 19(4):1757--1780, 2009.

\bibitem[NY17]{NY17}
Michael~J Neely and Hao Yu.
\newblock Online convex optimization with time-varying constraints.
\newblock {\em arXiv preprint arXiv:1702.04783}, 2017.

\bibitem[Roc70]{Roc70}
R.~Tyrrell Rockafellar.
\newblock {\em Convex Analysis}.
\newblock Princeton University Press, 1970.

\bibitem[Roc84]{Roc84}
R.~Tyrrell Rockafellar.
\newblock {\em Network flows and monotropic optimization}.
\newblock Wiley, New York, NY, 1984.

\bibitem[RW98]{RW98}
R.~Tyrrell Rockafellar and Roger J-B Wets.
\newblock {\em Variational Analysis}, volume 317.
\newblock Springer: Grundlehren der Math. Wissenschaften., Berlin, 1998.

\bibitem[Uza58]{Uza58}
H.~Uzawa.
\newblock Iterative methods in concave programming.
\newblock {\em Studies in Linear and Nonlinear Programming (K. Arrow, L.
  Hurwicz, and H. Uzawa, eds.)}, pages 154--165, 1958.

\bibitem[Van16]{Van16}
Lieven Vandenberghe.
\newblock {EE263C, UCLA}.
\newblock Lecture notes, 2016.

\bibitem[YNW17]{YNW17}
Hao Yu, Michael Neely, and Xiaohan Wei.
\newblock Online convex optimization with stochastic constraints.
\newblock In {\em Advances in Neural Information Processing Systems 30}, pages
  1428--1438, 2017.

\bibitem[Zin03]{Zin03}
Martin Zinkevich.
\newblock Online convex programming and generalized infinitesimal gradient
  ascent.
\newblock In {\em International Conference on Machine Learning}, 2003.

\end{thebibliography}
